\documentclass[12pt,oneside]{amsart}

\usepackage{amssymb, amsthm, enumerate, comment, url, color}
\usepackage{amsmath}

\newtheorem{theorem}{Theorem}
\newtheorem{lemma}[theorem]{Lemma}
\newtheorem{proposition}[theorem]{Proposition}
\newtheorem{corollary}[theorem]{Corollary}

\theoremstyle{definition}

\newtheoremstyle{remarks}%
{6pt}
{6pt}
{}
{}
{\bfseries}
{.}
{5pt}
{}%

\theoremstyle{remarks}
\newtheorem{remark}[theorem]{Remark}
\newtheorem{definition}[theorem]{Definition}

\numberwithin{theorem}{section}
\numberwithin{equation}{section}

\usepackage{geometry}

\newcommand{\lz}{\left(}
\newcommand{\pz}{\right)}
\renewcommand{\epsilon}{\varepsilon}
\renewcommand{\d}{\mathrm{d}}
\newcommand{\C}{\mathbb{C}}
\newcommand{\Z}{\mathbb{Z}}
\newcommand{\sumstar}{\sideset{}{^{*}}\sum}
\newcommand{\N}{\mathbb{N}}
\newcommand{\R}{\mathbb{R}}

\newcommand{\lab}{\left|}
\newcommand{\rab}{\right|}

\newcommand{\re}{\mathrm{Re}}
\newcommand{\im}{\mathrm{Im}}
\newcommand{\bfrac}[2]{\lz\frac{#1}{#2}\pz}

\renewcommand{\mod}[1]{\text{ (mod $#1$)}}
\newcommand{\odd}{\mathrm{\ odd}}

\renewcommand{\l}{\ell}
\renewcommand{\phi}{\varphi}
\newcommand{\leg}[2]{\left(\frac{#1}{#2}\right)}
\newcommand{\ch}{\mathrm{conv}}

\newcommand{\res}[1]{\underset{#1}{\mathrm{Res\ }}}

\begin{document}
	\title[First moment of quadratic Dirichlet $L$-functions with secondary terms]{First moment of quadratic Dirichlet $L$-functions with secondary terms}
	
	\author{M. \v Cech}
	
	\address[Martin \v Cech]
	{Charles University, Faculty of Mathematics and Physics, Department of Mathematical Analysis and Department of Algebra, Sokolovsk\'a 83, 18600 Praha 8, Czech Republic}
	\email{martin.cech@matfyz.cuni.cz}

	\begin{abstract}
		We study the first moment of primitive quadratic Dirichlet $L$-functions. Assuming the Riemann hypothesis and the generalized Lindel\"of hypothesis, we obtain an asymptotic formula at the central point with error $O(X^{1/4+\epsilon})$, and identify lower-order terms of size $X^{1/3}\log X$.
	\end{abstract}
\maketitle
\section{Introduction}

In this paper, we study the first moment of primitive quadratic Dirichlet $L$-functions, aiming to obtain the best possible error term. Our main result is the following:
\begin{corollary}\label{cor:main}
	Assume the Riemann hypothesis and the generalized Lindel\"of hypothesis and let $\Phi$ be a non-negative smooth weight supported in $[1/2,2]$ . Then for $X\geq 2$ and some (explicitly computed) linear polynomials $P_1(x), P_2(x)$, we have
	$$
	\sum_{\substack{d\geq 1\\ d\odd}} \mu^2(d)\Phi\bfrac dX L\lz1/2,\leg{8d}{\cdot}\pz = XP_1(\log X)+X^{1/3}P_2(\log X)+O(X^{1/4+\epsilon}).
	$$
\end{corollary} Corollary \ref{cor:main} is a direct consequence of the more general Theorem \ref{thm:moment} below. 

At the central point, the first moment was first computed by Jutila \cite{jutila}, who obtained an asymptotic formula (for the non-smooth moment) with error term $O(X^{3/4+\epsilon})$. Later, Goldfeld and Hoffstein \cite{goho} improved Jutila's error term, which is implicitly $O(X^{1/2+\epsilon})$ for the smooth moment. A similar error was also obtained by Young \cite{young_first_moment}, who used a different method whose main idea was to treat the square-free condition using a recursive argument. He used a similar technique to improve the error term in the third moment \cite{young_third_moment}, where he emphasized that the square-free condition is a substantial difficulty rather than a mere technicality, as is clear from comparing the sums
$$
\sum_{n\in\Z} f(n),\qquad \text{and} \qquad \sum_{n\in Z}\mu^2(n)f(n),
$$ where $f$ is a smooth Schwartz test function supported in $[-X,X]$. The first sum can be very accurately estimated using Poisson summation, while obtaining an error $O(X^{1/2-\delta})$ for some $\delta>0$ in the second sum is equivalent to the quasi-Riemann hypothesis.

Based on extensive computations, Alderson and Rubinstein \cite{alru} observed that the optimal error term seems to be $O(X^{1/4+\epsilon})$, but there appeared to be a bias in their computations.

Over function fields, where the analogue of the Riemann hypothesis is known, Florea \cite{florea_first_moment} managed to compute the first moment with (an analogue of) the best conjectured error $O(X^{1/4+\epsilon})$. She also found new, lower-order terms of size $X^{1/3}\log X$. Her treatment of the square-free condition in \cite[Lemma 2.2]{florea_first_moment} is based on a formula where polynomials are grouped together according to their degree, which does not seem to have a clear analogue over number fields.

Without the square-free condition, Gao and Zhao \cite{gazh_first_moment} (and implicitly Blomer \cite{blo}) under the Riemann hypothesis, computed the first moment at the central point with error $O(X^{1/4+\epsilon})$ by obtaining meromorphic continuation of the associated double Dirichlet series. The extra terms of size $X^{1/3}\log X$ do not appear there, so it is natural to ask whether a similar result can be obtained for primitive Dirichlet $L$-functions, and whether the extra terms appear. Both of these questions are affirmatively answered in Corollary \ref{cor:main}. 

For a better understanding where the lower-order terms come from and why they do not appear without the square-free condition, we refer to Remark \ref{rem:square-freeness}.

\smallskip

We now state our main results. We work with the family of characters $\chi_{8d}$ for odd square-free positive integers $d$, given by the Kronecker symbols $\leg{8d}{\cdot}$. A star on a sum $\sumstar$ indicates that it runs over square-free numbers. We denote the gamma factors in the functional equation of the corresponding $L$-functions by
$$
	X_+(s):=\bfrac{\pi}{8}^{s-1/2}\frac{\Gamma\bfrac{1-s}{2}}{\Gamma\bfrac s2}.
$$ Throughout the paper, $\Phi(x)$ denotes a non-negative smooth weight supported in $[1/2,2],$ and $\widetilde\Phi(s)$ denotes its Mellin transform.

\begin{theorem}\label{thm:moment}
	 Assume the generalized Lindel\"of hypothesis. Let $X\geq 2$, $\Phi$ be as above and define $$\beta:=\sup\{\re(\rho):\zeta(\rho)=0\}.$$
	\begin{enumerate}
		\item Fix $s\in\C$ with $1/3<\re(s)<1$, $s\neq 1/2$. Then for any $\epsilon>0,$
		\begin{equation}\label{eq:general moment result}
\begin{aligned}
				\sumstar_{d\odd}L(s,\chi_{8d})\Phi\bfrac dX&= X \widetilde\Phi(1)R_1(s)+ X^{\tfrac32-s}\widetilde\Phi(\tfrac32-s)X_+(s)R_1(1-s)\\
				&+X^{\tfrac12-\tfrac s3} \widetilde\Phi(\tfrac12-\tfrac s3)R_2(s)+X^{\frac{2-2s}{3}}\widetilde\Phi(\tfrac{2-2s}{3})X_+(s)R_2(1-s)\\
				&+O_\epsilon\lz X^{\frac{1-\re(s)}{2}+\epsilon}+X^{\frac\beta2+\epsilon}+X^{\frac{\beta+1}{2}-\re(s)+\epsilon}\pz,
\end{aligned}
		\end{equation} 
		where
		 $R_1(s)$  and $R_2(s)$ are explicit quantities defined in Theorem \ref{thm:MDS} and Proposition \ref{prop:residues}.
		\item We have for any $\epsilon>0$
		$$
			\sumstar_{d\odd} L(1/2,\chi_{8d})\Phi(d/X)=XP_1(\log X)+X^{1/3}P_2(\log X)+O_\epsilon(X^{\beta/2+\epsilon}),
		$$
		where $P_1(x), P_2(x)$ are linear polynomials given by
		$$
\begin{aligned}
				P_1(x)&=\frac{\widetilde\Phi(1)P\lz\tfrac 12\pz}{3\zeta(2)}\lz\frac{x}{2}+2\gamma+2\log 2+\frac{P'\lz\tfrac{1}{2}\pz}{P\lz\tfrac 12\pz}+\frac{\widetilde\Phi'(1)/\widetilde\Phi(1)-X_+'\lz\tfrac12\pz}{2}\pz, \\
				P_2(x)&= \widetilde\Phi\lz\tfrac13\pz T\lz\tfrac12\pz\lz \frac{x}{6}+2\gamma+\frac{\widetilde\Phi'\lz\tfrac13\pz}{6\widetilde\Phi\lz\tfrac13\pz}+\frac{T'\lz\tfrac12\pz}{T\lz\tfrac12\pz}-\frac{X_+'\lz\tfrac12\pz}{2}\pz,
\end{aligned}
		$$
		where $P(s)=\prod_{p>2}\lz1-\frac{1}{p^{2s}(p+1)}\pz$, and
		$$
			T(s) = Y(s)Q_2(s)Q(s),
		$$ where $Y$, $Q_2$, $Q$ are defined in Proposition \ref{prop:residues}.
	\end{enumerate}
\end{theorem}
As stated, our results are conditional under the generalized Lindel\"of hypothesis (GLH). In most places, it is enough to replace GLH by an estimate on average provided for example by Heath-Brown's quadratic large sieve \cite{HB_quadratic_sieve}, but the full GLH is required in one place of the proof of Proposition \ref{prop:continuation of B}, where we need to bound a product of arbitrarily many Dirichlet $L$-functions. It is possible to obtain a sligthly weaker unconditional result, which would in particular suffice to recover the error term $O(X^{1/2+\epsilon})$(see the discussion in Remark \ref{rem:no lindelof}).

The condition $\re(s)>1/3$ in part (1) is not necessary, it would be straightforward to extend the range of $\re(s)$. However, more (explicitly computable) main terms will arise in that case.

Part (2) of Theorem \ref{thm:moment} follows from part (1) by taking the limit $s\rightarrow 1/2$.

 To prove part (1), we use Perron's formula and investigate the resulting double Dirichlet series.
 
 Perron's formula gives
 \begin{equation}\label{eq:Perron}
 	\sumstar_{d\odd} L(s,\chi_{8d})\Phi(d/X)=\frac{1}{2\pi i}\int_{(c)}A(s,w)\tilde \Phi(s) X^sds,
 \end{equation} where
 $$
 	A(s,w)=\sumstar_{d\odd}\frac{L(s,\chi_{8d})}{d^w}.
 $$
 
 Theorem \ref{thm:moment} (1) is a consequence of the following result.

\begin{theorem}\label{thm:MDS}
	Assume the generalized Lindel\"of hypothesis. Then the double Dirichlet series $A(s,w)$ has a meromorphic continuation to the region
	$$
		\{\re(s+w)>1/2,\ \re(s+2w)>1,\ \re(w)>0\}.
	$$
The only possible poles of $\zeta(2w)\zeta(2s+2w-1)A(s,w)$ in this region occur at $w=1$, $s+w=3/2$, $s+(2j+1)w=3/2$ and $2js+(2j+1)w=j+1$ for $j\in \N$. If $0<\re(s)<3/2$ and $s\neq 1/2$, it has the following residues:
$$
\begin{aligned}
		\res{w=1}A(s,w)&= R_1(s),\\
	\res{w=3/2-s}A(s,w)&= X_+(s) R_1(1-s), \\
	\res{w=1/2-s/3}A(s,w)&=R_2(s), \\
	\res{w=2/3-2s/3}A(s,w)&=X_+(s)R_2(1-s),
\end{aligned}
$$ where $R_1(s), R_2(s)$ are explicitly computed in Proposition \ref{prop:residues}.
Moreover, $A(s,w)$ is polynomially bounded in vertical strips in this region.
\end{theorem}
Polynomial boundedness in vertical strips in our setting is defined in definitions \ref{def:bounded} and \ref{def:bounded meromorphic}.

\begin{proof}[Proof of Theorem \ref{thm:moment} assuming that Theorem \ref{thm:MDS} holds] Let $\beta$ be as in Theorem \ref{thm:moment}.
	
	Proof of part (1): assume that $s$ is as in the Theorem. By Perron's formula, we have 
	$$
		\sumstar_{d\geq 1}L(s,\chi_{8d})\Phi(d/X)=\frac{1}{2\pi i}\int_{(2)}A(s,w)\widetilde\Phi(w)X^w\d w.
	$$
	Let $$c=\max\{\beta/2, (\beta+1)/2-\re(s),1/2-\re(s)/2\},$$ and shift the integral to the line $\re(w)=c+\epsilon$. Theorem \ref{thm:MDS} assures that we are within the region of meromorphic continuation of the integrand. Thanks to the smoothness of $\Phi$ and polynomial boundedness in vertical strips of $A(s,w),$ the shifted integral is absolutely convergent and bounded by $O(X^{c+\epsilon})$, which gives the claimed error term. The main terms come from the residues at the poles at $w=1$, $w=3/2-s$, $w=1/2-s/3$ and $w=2/3-2s/3$ (note that some of them were not necessarilly crossed during the integral shift, but the formula remains true anyway); we now argue that no other poles are crossed. 
	
	Since $c>\max\{\beta/2,(\beta+1)/2-\re(s)\}$, there is no contribution from the possible poles at $w=\rho/2$ or $w=(\rho+1)/2-s$.
	
	We also did not cross any of the poles at $w=\frac{3}{2(2j+1)}-\frac{s}{2j+1}$ for $j\geq 2$. This is clear if $j\geq 3$ since in that case the pole is on the line $\re(w)<\frac{3}{14}<\frac14\leq \frac{\beta}{2}<c+\epsilon$. If $j=2$, since $\re(s)>1/3$, we also find that this pole satisfies $\re(w)=\frac{3}{10}-\frac{\re(s)}{5}<\frac14.$
	
	It remains to show that we did not cross any of the possible poles at $w=\frac{j+1}{2j+1}-\frac{2js}{2j+1}$ for $j\geq 2$. These poles are at the points $w=\frac12+\frac{1-4js}{4j+2}=\frac{1}{2}-\frac{s}{2}+\frac{(1-2j)s+1}{4j+2},$ but the assumptions $\re(s)>1/3$ and $j\geq 2$ imply that $(1-2j)\re(s)+1<0.$ Therefore these poles satisfy $\re(w)<\frac12-\frac{\re(s)}{2}\leq c$, so they were not crossed during the integral shift.
	
	\smallskip
	
	Proof of part (2): this part follows from taking the limit as $s\rightarrow 1/2$ in (1). Note that $R_1(s), R_2(s)$ have a simple pole at $s=1/2$, but we find that the poles of the first two and second two terms in \eqref{eq:general moment result} cancel.
	
	To find the result, we expand each function in the main term as a Laurent series centered at $s=1/2$. Writing $R_1(s)=\frac{2}{3\zeta(2)}\zeta(2s)P(s)\lz1-\frac{1}{2^{2s}}\pz$ as in Proposition \ref{prop:residues}, we find the following Laurent expansions.
	$$
		\begin{aligned}
			X\widetilde\Phi(1)R_1(s)&=\frac{X\widetilde\Phi(1)2}{3\zeta(2)}\lz\frac{1}{2s-1}+\gamma+\dots\pz\lz P\lz\tfrac12\pz+(s-\tfrac12)P'\lz\tfrac12\pz+\dots\pz\\
			&\cdot\lz\frac12+(s-\tfrac12)\log 2+\dots\pz\\
			&=\frac{X\widetilde\Phi(1)P\lz\tfrac12\pz}{3\zeta(2)}\lz\frac{1}{2s-1}+\gamma +\frac{P'\lz\tfrac12\pz}{2P\lz\tfrac12\pz}+\log 2+O\lz s-\tfrac12\pz\pz,
		\end{aligned}
	$$ and
	$$
\begin{aligned}
	&X^{\frac32-s}\widetilde\Phi(\tfrac32-s)X_+(s)R_1(1-s)=\frac{2X}{3\zeta(2)}\lz\widetilde\Phi(1)-(s-\tfrac12)\widetilde\Phi'(1)+\dots\pz\\
	&\cdot \lz1-(s-\tfrac12)\log X+\dots\pz\lz1+(s-\tfrac12)X_+'\lz\tfrac12\pz+\dots\pz\lz P\lz\tfrac12\pz-(s-\tfrac12)P'\lz\tfrac12\pz+\dots\pz\\
	&\cdot \lz\frac{1}{1-2s}+\gamma+\dots\pz\lz\tfrac12-\log 2( s-\tfrac12)+\dots\pz\\
	&=\frac{X\widetilde\Phi(1) P\lz\tfrac12\pz }{3\zeta(2)}\lz\frac{1}{1-2s}+\gamma +\tfrac{ P'\lz\tfrac12\pz}{2P\lz\tfrac12\pz}+\tfrac{\log X -X_+'\lz\tfrac12\pz}{2}+\tfrac{\widetilde\Phi'(1)}{2\widetilde\Phi(1)}+\log 2+O\lz s-\tfrac12\pz\pz
\end{aligned}
	$$
	Taking these terms together in the limit as $s\rightarrow 1/2$ leads to the main term $XP_1(\log X)$. 
	
	The second term can be obtained similarly from the following Laurent expansions, where we write $R_2(s)=\zeta(2s)T(s)$, with $T(s)=Y(s)Q_2(s)Q(s)$ with the notation from Proposition \ref{prop:residues}:
	$$
		\begin{aligned}
			X^{\tfrac12-\tfrac s3}&\widetilde\Phi\lz\tfrac12-\tfrac s3\pz R_2(s)=X^{\tfrac13}\lz1+\lz\tfrac16-\tfrac s3\pz\log X+\dots\pz\lz \widetilde\Phi(\tfrac13)+\lz\tfrac16-\tfrac s3\pz\widetilde\Phi'\lz\tfrac 13\pz+\dots\pz\\
			&\qquad\qquad\qquad\quad\cdot\quad\lz \frac{1}{2s-1}+\gamma+\dots\pz \lz T\lz\tfrac12\pz +(s-\tfrac12)T'\lz\tfrac12\pz+\dots\pz\\
			&=X^{\frac13} \widetilde\Phi \lz\tfrac 13\pz T\lz\tfrac 12\pz\lz \frac{1}{2s-1}+\gamma-\frac{\log X}{6}-\frac{\widetilde\Phi'\lz\tfrac 13\pz}{6\widetilde\Phi\lz\tfrac 13\pz}+\frac{T'\lz\tfrac 12\pz}{2T\lz\tfrac 12\pz}+O\lz s-\tfrac12\pz\pz,
		\end{aligned}
	$$ and
	$$
		\begin{aligned}
			&X^{\tfrac23-\tfrac {2s}{3}}\widetilde\Phi\lz\tfrac 23-\tfrac{2s}{3}\pz X_+(s) R_2(1-s)=X^{\tfrac13}\lz 1+\lz\tfrac 13-\tfrac{2s}{3}\pz\log X+\dots\pz \\
			&\cdot\lz \widetilde\Phi\lz\tfrac13\pz+\lz\tfrac13-\tfrac{2s}{3}\pz\widetilde\Phi'\lz\tfrac 13\pz+\dots\pz \lz 1+\lz s-\tfrac12\pz X_+'\lz\tfrac 12\pz+\dots\pz\\
			&\cdot \lz\frac{1}{1-2s}+\gamma+\dots\pz \lz T\lz\tfrac 12\pz-(s-\tfrac12)T'\lz\tfrac12\pz+\dots\pz \\
			&=X^{\frac13} \widetilde\Phi\lz\tfrac13\pz T\lz\tfrac12\pz \lz \frac{1}{1-2s}+\tfrac{\log X}{3}+\tfrac{\widetilde\Phi'\lz\tfrac13\pz}{3\widetilde\Phi \lz\tfrac13\pz}-\tfrac{X_+'\lz\tfrac12\pz}{2}+\gamma+\tfrac{T'\lz\tfrac12\pz}{2T\lz\tfrac12\pz}+O\lz s-\tfrac12\pz\pz.
		\end{aligned}
	$$
	Taking these two terms together in the limit as $s\rightarrow 1/2$ gives the second main term.
\end{proof}

Our proof of Theorem \ref{thm:MDS} uses the usual strategy outlined for example by Diaconu, Goldfeld and Hoffstein in \cite{dgh}, where one uses a region of absolute convergence, functional equations of the multiple Dirichlet series (MDS) and Bochner's tube Theorem (see Theorem \ref{thm:bochner}) to obtain a meromorphic continuation of the MDS in question. There are two functional equations of $A(s,w)$: a straightforward one uses the functional equation of $L(s,\chi)$ in the definition, and a less clear one comes after expanding each $L$-function into a Dirichlet series and exchanging sums. This new expression will contain the Dirichlet $L$-functions $L(w,\chi)$ and one would hope to use the functional equation of these, but an issue is that the characters $\chi$ appearing during this procedure are not always primitive. Thus one only has a heuristic functional equation of this form and the question is whether it can be obtained rigorously. A usual strategy is to modify the MDS in such a way that the heuristic functional equation can be proved rigorously.

We use an approach introduced in \cite{ratios}, where we applied a functional equation of non-primitive $L$-functions that related them to a Dirichlet series with corresponding Gauss sums (see Proposition \ref{prop:general fe}).

This approach was also used by the author in \cite{moments} to obtain a meromorphic continuation of MDS corresponding to the general $k$-th moment of real Dirichlet $L$-functions, with a the hope that a more careful treatment of the cases with $k=1,2,3$ would lead to a better region of meromorphic continuation and thus a better error term in the asymptotic formula for the moment. In this paper, we treat the case $k=1$, and it is plausible that a similar improvement would be possible for the cases $k=2$ and $k=3$. The expected best possible error terms using these techniques would be $O(X^{1/2+\epsilon})$ for the second moment and $O(X^{3/4+\epsilon})$ for the third moment, which would match that of Young \cite{young_third_moment}. In particular, one should not hope to obtain a better error term than $O(X^{3/4+\epsilon})$ for $k=3$ due to the presence of lower order terms of size $X^{3/4}$, which come from a very different source than our lower-order terms (see also the discussion below Corollary 1.5 in \cite{moments}).

\smallskip

After stating some preliminaries in Section \ref{sec:prelim}, we will prove the basic properties of $A(s,w)$ in \ref{sec:basic}. In Section \ref{sec:fe in w}, we will then focus on obtaining a bigger region of meromorphic continuation by using the extra functional equation. Finally, to obtain explicit polynomials $P_1,P_2$ in the main results, we explicitly compute the residues of $A(s,w)$ in the required poles in Section \ref{sec:residues}; these computations are somewhat tedious, but straightforward.

\section*{Acknowledgments}
The author was supported by the Charles University grants PRIMUS/24/SCI/010 and PRIMUS/25/SCI/017.

\section{Preliminaries}\label{sec:prelim}
Throughout the paper, $\epsilon$ will denote an arbitrarily small positive real number that can be different at each occurrence. All implied constants are allowed to depend on $\epsilon$.

We define the arithmetic function $a_t(n)$ by \begin{equation}\label{eq:def of a_t(n)}
	a_t(n)=\prod\limits_{p|n}\lz 1-\frac{1}{p^t}\pz^{-1},
\end{equation} and note that if $\re(t)>0$, then
\begin{equation}\label{eq:bound for a_t(n)}
	|a_{t}(n)|\ll_t n^\epsilon.
\end{equation}

\subsection{Dirichlet characters and Gauss sums}
For $n\in\Z$, we will denote by $\chi_n$ the Kronecker symbol $\leg{n}{\cdot}$. It is not always a Dirichlet character, but if $d$ is an odd and squarefree positive integer, then $\chi_{8d}$ is an even primitive character of conductor $8d$. Moreover, $\chi_{4n}$ is an even Dirichlet character modulo $4n$ for any positive integer $n$. 

We will also work with the Kronecker symbol $\leg{\cdot}{n}$ for a positive integers $n$. This is always a Dirichlet character, whose modulus is $n$ if $n\equiv 0,1,3\mod 4$, and $4n$ if $n\equiv 2\mod 4$. If $n$ is odd, then the characters $\leg{\cdot}{n}, \leg{\cdot}{4n}$ are even if $n\equiv 1\mod 4$ and odd if $n\equiv -1\mod 4$. If $n$ is squarefree, then $\leg{\cdot}{n}$ is primitive.

By $\psi_1,\psi_{-1}$, we denote respectively the principal and non-principal character modulo $4.$ For $j=\pm 2$, we also define $\psi_{j}:=\leg{4j}{\cdot}$, which are both primitive characters modulo $8$. Finally, by $\psi_0$ we denote the ``primitive character modulo 1'', which is the constant function 1.

\smallskip

For a Dirichlet character $\chi$ modulo $n$, the Gauss sum $\tau(\chi,k)$ is defined by
$$
\tau(\chi,k)=\sum_{j\mod n}\chi(j)e\bfrac{jk}{n},
$$ where we used the standard notation $e(x):=e^{2\pi ix}.$ For $\chi$ modulo $n$ and integers $k_1,k_2$ with $(n,k_1)=1$, we have $\tau(\chi,k_1k_2)=\overline\chi(k_1)\tau(\chi,k_2)$:
\begin{equation}\label{eq:gauss sums second coord}
	\sum_{j\mod n}\chi(j)e\bfrac{jk_1k_2}{n}=\sum_{j\mod n}\chi(j\overline{k_1})e\bfrac{jk_2}{n}=\overline\chi(k_1)\tau(\chi,k_2).
\end{equation}

A direct computation using the Chinese remainder Theorem also gives that for odd $n$ and any integer $k$, we have
\begin{equation}\label{eq:Gauss sums with 4}
	\tau\lz\leg{\cdot}{4n},k\pz=\tau\lz\psi_1,k\pz\tau\lz\leg{\cdot}{n},k\pz=\begin{cases}
		0,&\hbox{$k$ odd,}\\
		-2\tau\lz\leg{\cdot}{n},k\pz,&\hbox{$k\equiv 2\mod 4,$}\\
		2\tau\lz\leg{\cdot}{n},k\pz,&\hbox{$k\equiv 0\mod 4$.}
	\end{cases}	
\end{equation}

\smallskip

For the Jacobi symbol $\leg{\cdot}{n}$, we further define the modified Gauss sums by
$$
	\begin{aligned}
		G\lz\leg{\cdot}{n},k\pz&=\lz\frac{1-i}{2}+\leg{-1}{n}\frac{1+i}{2}\pz\tau\lz\leg{\cdot}{n},k\pz\\
		&=\begin{cases}
			\tau\lz\bfrac{\cdot}{n},k\pz,&\hbox{if $n\equiv 1\mod 4$,}\\
			-i\tau\lz\bfrac{\cdot}{n},k\pz,&\hbox{if $n\equiv 3\mod 4$,}
		\end{cases}
	\end{aligned}
$$ or equivalently
\begin{equation}\label{eq:tau as G}
	\tau\lz\leg{\cdot}{n},k\pz=\begin{cases}
		G\lz\leg{\cdot}{n},k\pz,&\hbox{if $n\equiv 1\mod 4$,}\\
		iG\lz\leg{\cdot}{n},k\pz,&\hbox{if $n\equiv3\mod 4$}.
	\end{cases}
\end{equation}
The advantage is that $G\lz\leg{\cdot}{n},k\pz$ are multiplicative in the $n$-variable: by \cite[Lemma 2.3]{sound} for $(m,n)=1$, we have
$$
	G\lz\leg{\cdot}{m},k\pz G\lz\leg{\cdot}{n},k\pz=G\lz\leg{\cdot}{mn},k\pz,
$$
and we also have the following explicit evaluation: if $p$ is an odd prime and $p^a||k$, then
\begin{equation}\label{eq:Gauss sum evaluation}
	G\lz\leg{\cdot}{p^j},k\pz=\begin{cases}
		\phi(p^j),&\hbox{if $j\leq a$, $j$ even,}\\		0,&\hbox{if $j\leq a$, $j$ odd,}\\
		-p^a,&\hbox{if $j=a+1$, $j$ even,}\\
		\leg{k p^{-a}}{p}p^a\sqrt p,&\hbox{if $j=a+1$, $j$ odd,}\\
		0,&\hbox{if $j\geq a+2$.}\\
	\end{cases}
\end{equation}

\subsection{Dirichlet $L$-functions} Throughout the paper, $\zeta(s)$ and $L(s,\chi)$ denote respectively the Riemann zeta function and the Dirichlet $L$-function. A subscript $\zeta_{(k)}$ or $L_{(k)}$ denotes the corresponding $L$-function with Euler factors for $p|k$ removed.

If $\chi$ is an even primitive Dirichlet character of conductor $n$, then we have the functional equation
\begin{equation}\label{eq:fe}
	L(s,\chi)=n^{1/2-s} Y_+(s)L(1-s,\overline\chi),\qquad \text{where} \qquad Y_+(s)=\pi^{s-1/2}\frac{\Gamma\bfrac{1-s}{2}}{\Gamma\bfrac s2}.
\end{equation}

Our main result is conditional on the Lindel\"of hypothesis for Dirichlet $L$-functions, which is the following estimate: for any primitive Dirichlet character $\chi$ of conductor $n$ and arbitrary $\sigma\geq 1/2$, $t\in\R$, we have
\begin{equation}\label{eq:LH}
	|L(\sigma+it,\chi)|\ll_{t}n^{\epsilon}.
\end{equation}
 This implies a similar bound for non-primitive characters by writing the corresponding $L$-function in terms of the primitive one.

In other instances, we can replace \eqref{eq:LH} with the corresponding bound on average, which is provided by a large sieve estimate due to Heath-Brown \cite{HB_quadratic_sieve}:
\begin{lemma}\label{le:lindelof on average}
	\begin{enumerate}
		\item For any $\sigma,t\in \R$, we have 
		$$
		\sumstar_{d\leq X}|L(\sigma+it,\chi_{8d})|\ll_{t,\sigma} X^{\max\{1,3/2-\sigma\}+\epsilon}.
		$$
		\item For any $\sigma\geq 0$ and $t\in\R$, we have 
		$$
			\sum_{n\leq X}\lab L\lz\sigma+it,\leg{\cdot}{n}\pz\rab\ll_{t,\sigma} X^{\max\{1,3/2-\sigma\}+\epsilon}.
		$$
	\end{enumerate}
\end{lemma}
\begin{proof}
	(a): by \cite[Theorem 2]{HB_quadratic_sieve}, and Holder's inequality, we have for any $\sigma\geq 1/2$
	$$
		\sumstar_{d\leq X}|L(\sigma+it,\chi_d)|\ll_t X^{1+\epsilon}. 
	$$
	By \eqref{eq:fe}, if $\sigma<1/2$ we have
	$$
		\sumstar_{d\leq X}|L(\sigma+it,\chi_d)|=|Y_+(s)|\sumstar_{d\leq X}|L(1-\sigma-it,\chi_d)|d^{1/2-\sigma}\ll_{t,\sigma} X^{3/2-\sigma+\epsilon}. 
	$$
	Altogether, we obtain
	$$
		\sumstar_{d\leq X}|L(\sigma+it,\chi_d)|\ll_t X^{\max\{1,3/2-\sigma\}+\epsilon},
	$$ which finishes the proof of part (a).
	
	(b): we write $n=n_0n_1^2$ with $n_0$ square-free, so for any $s\in\C$ with $\re(s)\geq 0$, we have
	$$
		\lab L\lz s,\leg{\cdot}{n}\pz\rab=\lab L\lz s,\leg{\cdot}{n_0}\pz\rab\lab\prod_{p|2n_1}\lz1-\frac{\chi_{n_0}(p)}{p^s}\pz\rab\ll |L(s,\chi_{n_0})| d(n_1).
	$$ Note that part (a) holds with the same proof if the family of characters is replaced by an arbitrary set of primitive quadratic characters of conductor $\ll X$ because \cite[Theorem 2]{HB_quadratic_sieve} holds for the sum of all such characters. Therefore if $\sigma\geq0$,
	$$
		\begin{aligned}
			\sum_{n\leq X}\lab L\lz\sigma+it,\leg{\cdot}{n}\pz\rab&\ll \sumstar_{n_0\leq X}\lab L\lz\sigma+it,\leg{\cdot}{n_0}\rab\pz\sum_{n_1^2\leq X/n_0}d(n_1)\\
			&\ll X^{1/2+\epsilon}\sumstar_{n_0\leq X}\frac{\lab L\lz\sigma+it,\leg{\cdot}{n_0}\pz\rab}{n_0^{1/2}}\\
			&\ll X^{\max\{1,3/2-\sigma\}+\epsilon},
		\end{aligned}
	$$ where the last estimate follows from partial summation.
\end{proof} As remarked during the proof of part (b), analogues of the lemma above hold with the same proof when we sum over other families of real characters with conductor $\ll X$.

\smallskip

We will make crucial use of the following functional equation for non-primitive characters, which is \cite[Proposition 2.3]{ratios} (with a slightly different normalization).

\begin{proposition}\label{prop:general fe}
Let $\chi$ be an arbitrary character modulo $n$. Then
\begin{equation}\label{eq:general fe}
	L(s,\chi)=n^{1/2-s}Y_{\pm}(s)K(1-s,\chi),
\end{equation} where
$$
	\begin{aligned}
		K(s,\chi)&=\sum_{k=1}^\infty\frac{\tau(\chi,k)/\sqrt n}{k^s},\\
		Y_{\pm}(s)&=\begin{cases}
			Y_+(s):=\pi^{s-1/2}\frac{\Gamma\bfrac{1-s}{2}}{\Gamma\bfrac s2},&\hbox{if $\chi$ is even,}\\
			Y_-(s):=-i\pi^{s-1/2}\frac{\Gamma\bfrac{2-s}{2}}{\Gamma\bfrac {1+s}2},&\hbox{if $\chi$ is odd.}
		\end{cases}
	\end{aligned}
$$
\end{proposition} 
From this functional equation, we find that $K(s,\chi)$ has very similar properties as $L(s,\chi)$. Mainly, it has a meromorphic continuation to all $s\in \C$ with the only pole at $s=1$ if $\chi$ is a principal character. Moreover, the poles of $Y_{\pm}(s)$ get cancelled by the ``trivial zeros'' of $K(s,\chi)$.

We also have the following Lindel\"of bound on average.
\begin{lemma}\label{le:lindelof for K}
	For any $\sigma>0$ and $t\in\R$, we have
$$
		\sum_{n\leq X}\lab K\lz\sigma+it,\leg{\cdot}{n}\pz\rab\ll_{t,\sigma} X^{{\max \{1,3/2-\sigma\}}+\epsilon}.
	$$
\end{lemma}
\begin{proof}
	If $\sigma>1$, then $K(\sigma+it,\chi)$ is absolutely convergent for any character $\chi\mod n$ since $|\tau(\chi,k)|\ll_n 1$, and the bound follows.
	
	For $0<\sigma\leq 1$, Proposition \ref{prop:general fe} gives $ K\lz\sigma+it,\leg{\cdot}{n}\pz=n^{1/2-s}Y_{\pm}(s)L\lz1-\sigma-it,\leg{\cdot}{n}\pz$ and $1-\sigma\geq 0$, so by Lemma \ref{le:lindelof on average} and partial summation,
	$$
		\begin{aligned}
			\sum_{n\leq X}\lab K\lz\sigma+it,\leg{\cdot}{n}\pz\rab&=\sum_{n\leq X}\frac{\lab L\lz 1-\sigma-it,\leg{\cdot}{n}\pz\rab}{n^{\sigma-1/2}}|Y_{\pm}(s)|\\
			&\ll_{t,\sigma}X^{\max\{1,1/2+\sigma\}+\epsilon}X^{1/2-\sigma}\\
			&=X^{\max\{1,3/2-\sigma\}+\epsilon}.
		\end{aligned}
	$$
\end{proof}

\subsection{multivariable complex analysis} 

In our application of Perron's formula in \eqref{eq:Perron}, we can consider $s$ as fixed and integrate only in the $w$-variable, thus thinking of $A(s,w)$ as a function of one complex variable $w$. However, the strength of our result mainly depends on our ability to obtain a meromorphic continuation of the function $A(s,w)$, for which it is useful to view it as a function of two complex variables.

 In general, a function of several complex variables is holomorphic if and only if it is holomorphic in each variable separately, so most of the basic theorems from single-variable complex analysis have their multivariable analogues. Without further mention, we will make frequent use of the uniqueness of the holomorphic (or meromorphic) continuation by performing computations with Dirichlet series in regions of absolute convergence and then extending the formulas to bigger regions, where all quantities have a holomorphic (or meromorphic) extension.
 
 One difference between functions of one or several complex variables is in the structure of their singularities. The singularities of functions of several complex variables are not isolated, but lie on some hypersurface (they are the zero locus of a holomorphic function), which can be very complicated. However, in our case, the singularities are simple enough -- they will be given by lines. Moreover, as mentioned above, the integration in \eqref{eq:Perron} is only over one variable, so we can think of $s$ fixed and use the usual  single-variable residue theorem to compute the integral.
 
 Another important difference is that in more than one variable, we sometimes obtain holomorphic continuation for free (see Theorem \ref{thm:bochner}).

A general reference for the theory of multivariable complex analysis is \cite{hor}.

\begin{definition}
	An open set $R\subset \C^n$ is a \emph{domain of holomorphy} if there are no open sets $R_1,R_2\subset\C^n$ such that $\emptyset\neq R_1\subset R\cap R_2,$ $R_2$ is connected and not contained in $R$, and for any holomorphic function $f$ on $R$, there is a function $f_2$ holomorphic on $R_2$ such that $f=f_2$ on $R_1.$ 
\end{definition}
Open balls $B(c,r)$ centered in $c$ of radius $r$ are domains of holomorphy.

The following is a generalization of vertical strips in $\C^n$.
\begin{definition}
	A set $T\subset\C^n$ is a \emph{tube domain} if there is a connected  set $U\subset\R^n$ such that $T=U+i\R^n=\{z\in\C^n:\ \re(z)\in U\}.$
	The set $U$ is called the \emph{base} of $T$.
	
	We say that a tube domain is \emph{open, bounded, compact}, etc. if its base has this property. 
	
	A \emph{tube subdomain} of $T$ is a set $S\subset T$  that is also a tube domain.
\end{definition}

The following is (a generalization of) Bochner's Tube Theorem \cite{Boc}.
\begin{theorem}\label{thm:bochner}
	An open tube domain is a domain of holomorphy if and only if it is convex.	
\end{theorem}
We denote the convex hull of $T$ by $\ch(T)$. In particular, every holomorphic function on an open tube domain $T$ has a holomorphic continuation to $\ch (T)$. The same result holds for functions on $T$ that are meromorphic if there are finitely many functions $g$ holomorphic on $T$ such that the product $f g$ is holomorphic on $T$ (note that for a general meromorphic function, this property holds only locally in $T$).

\smallskip

Now we collect some results that will enable us to control the size of the meromorphic continuation provided by Bochner's Tube Theorem.

\begin{theorem}\label{Preimage is domain of holomorphy}
	Let $R_1\subset\C^m,R_2\subset\C^n$ be domains of holomorphy, and $f~:~R_1~\rightarrow~\C^n$ a holomorphic map. Then $$
		R=f^{-1}(R_2)=\{z\in R_1:f(z)\in R_2\}
	$$ is a domain of holomorphy.
\end{theorem}

\begin{proposition}\label{prop:extending inequalities}
	Assume that $T\subset \C^n$ is an open tube domain, $g,h:T\rightarrow \C$ are holomorphic functions, and let $\tilde g,\tilde h$ be their holomorphic continuation to $\ch(T)$. If  $\lvert g(z)\rvert \leq \lvert h(z)\rvert $ for all $z\in T$, and $h(z)$ is nonzero in $T$, then also $\lvert \tilde g(z)\rvert \leq \lvert \tilde h(z)\rvert $ for all $z\in \ch (T)$.
\end{proposition}
\begin{proof}
	Since $h(z)$ is nonzero, the function $f(z)=g(z)/h(z)$ is holomorphic in $T$, and hence has a holomorphic continuation $\tilde f$ to $\ch (T)$. By our assumptions, $\lvert \tilde f(z)\rvert \leq 1$ for all $z\in T,$ so $T\subset \tilde f^{-1}(B(0,1))$. However, by Theorem \ref{Preimage is domain of holomorphy}, $\tilde f^{-1}(B(0,1))$ is a domain of holomorphy, so it is the whole $\ch (T).$
\end{proof}

We will only work with functions of two complex variables, so we restrict to these for now.
\begin{definition}\label{def:bounded}
	Let $T\subset \C^2$ be a tube domain and $f(s,w):T\rightarrow\C$ be a holomorphic function on $T$. We say that $f$ is \emph{bounded in vertical strips in $T$} if for every compact tube subdomain $S\subset T$, there are constants $\alpha,\beta$ such that for all $(s,w)\in S$, we have
	$$
	|f(s,w)|\ll (1+|s|)^{\alpha} (1+|w|)^{\beta}.
	$$
\end{definition}

We will also need a similar definition for meromorphic functions. In general, the pole structure of meromorphic functions in several complex variables can be very complicated, since they are given by zero sets of holomorphic functions. However, we can restrict our attention to those whose poles are lines. We thus make the following definition.

\begin{definition}\label{def:bounded meromorphic}
	Let $T\subset \C^2$ be a tube domain and $f:T\rightarrow\C$ be a meromorphic function on $T$. We say that $f$ is \emph{bounded in vertical strips in $T$} if for every compact subdomain $S\subset T$, there are finitely many linear functions $\l_1(s,w),\dots,\l_n(s,w)$ such that the function $f(s,w)\l_1(s,w)\dots \l_n(s,w)$ is holomorphic and bounded in vertical strips in $S$.
\end{definition}

We are now ready to show that if $f(s,w)$ is meromorphic and bounded in vertical strips on a tube domain $T$, then its extension to $\ch(T)$ is also bounded in vertical strips. We do so by transferring to the single-variable case and using the Phragm\'en-Lindel\"of principle there.

First, we prove an auxiliary lemma showing that the extension is of finite order.

\begin{lemma}\label{le:first order in convex hulls}
	Let $T_1,T_2\subset\C^2$ be bounded open tube domains such that $T_1\cap T_2\neq\emptyset$, and assume that for each $\epsilon>0$,
	$$
	\begin{aligned}
		|g(s_1,w_1)|&\ll_\epsilon e^{\epsilon (|\im(s_1)|+|\im(w_1)|)}, \quad (s_1,w_1)\in T_1,\\
		|g(s_2,w_2)|&\ll_\epsilon e^{\epsilon (|\im(s_2)|+|\im(w_2)|)}, \quad (s_2,w_2)\in T_2.
	\end{aligned}
	$$
	Then also for every $\epsilon>0$,
	$$
	|g(s,w)|\ll_\epsilon e^{\epsilon(|\im(s)|+|\im(w)|)},\quad (s,w)\in T:=\ch(T_1,T_2).
	$$
\end{lemma}
\begin{proof}
	Let $B_1,B_2\subset \R^2$ denote the bases of $T_1,T_2$ respectively. Since they are bounded, let $K\in \R$ be such that whenever $(a,b)\in B_1\cup B_2$, then $|a|, |b|, |a+b|<K$. Let now $0<\epsilon<1/K$ be arbitrary, and consider the function $$\begin{aligned}
		h_\epsilon(s,w)&:=e^{\epsilon i(s+w)}+e^{\epsilon i (s-w)}+e^{\epsilon i (w-s)}+e^{\epsilon i(-s-w)}\\
		&=\lz e^{\epsilon is}+e^{-\epsilon is}\pz\lz e^{\epsilon iw}+e^{-\epsilon iw}\pz.
	\end{aligned}$$
Note that $h_\epsilon(s,w)$ is holomorphic in $\C^2$, and $h_\epsilon(s,w)=0$ if and only if $\epsilon s=\pi/2+k\pi$ or $\epsilon w=\pi/2+k\pi$ for some $k\in \Z$. By our choice of $\epsilon$, for $(s,w)\in T_1\cup T_2$ we have $\epsilon|\re(s)|, \epsilon|\re(w)|<1<\pi/2$, so $h_\epsilon(s,w)\neq 0$ for $(s,w)\in T_1\cup T_2$. Moreover, $|h_\epsilon (s,w)|\asymp_\epsilon e^{\epsilon(|\im(s)|+|\im(w)|)}$ in $T_1\cup T_2$, so we have $|g(s,w)|\ll_\epsilon |h_\epsilon(s,w)|$ on $T_1\cup T_2$. The conclusion now follows from Proposition \ref{prop:extending inequalities}.
\end{proof}

\begin{lemma}\label{le:bounded in vertical strips holom}
	Let $R_1,R_2\subset \C^2$ be open tube domains such that $R_1\cap R_2\neq \emptyset$ and $g(s,w)$ be a function that is holomorphic on $R:=\ch(R_1,R_2).$ Assume that there are $\alpha_1,\beta_1,\alpha_2,\beta_2>0$ such that
	\begin{equation}\label{eq:bounds in vertical strips assumptions}
		\begin{aligned}
		|g(s_1,w_1)|&\ll (1+|\im(s_1)|)^{\alpha_1}(1+|\im(w_1)|)^{\beta_1}, \quad (s_1,w_1)\in R_1,\\
		|g(s_2,w_2)|&\ll (1+|\im(s_2)|)^{\alpha_2}(1+|\im(w_2)|)^{\beta_2}, \quad (s_2,w_2)\in R_2.
	\end{aligned}
	\end{equation}
	Let $(s,w)\in R$ be arbitrary and let $(s_1,w_1)\in R_1$, $(s_2,w_2)\in R_2$ and $0\leq t\leq 1$ be such that $(s,w)=t(s_1,w_1)+(1-t)(s_2,w_2)$. Then
\begin{equation}\label{eq:bound in vertical strips final}
	\begin{aligned}
		|g(s,w)|&\ll (1+|\im (s)|)^{t\alpha_1+(1-t)\alpha_2}(1+|\im(w)|)^{t\beta_1+(1-t)\beta_2}\\
		&\cdot (1+|\re(s_1)-\re(s_2)|)^{t\alpha_1+(1-t)\alpha_2}(1+|\re(w_1)-\re(w_2)|)^{t\beta_1+(1-t)\beta_2},
\end{aligned}
\end{equation}
	 where the implied constant depends only on the implied constants in \eqref{eq:bounds in vertical strips assumptions}.
	
	In particular, if $g(s,w)$ is bounded in vertical strips in $R_1$ and $R_2$, then it is also bounded in vertical strips in $R$.
\end{lemma}

\begin{proof}
	Let $(s,w)\in R$ and find $(s_1,w_1)\in R_1$, $(s_2,w_2)\in R_2$ such that $(s,w)=t(s_1,w_1)+(1-t)(s_2,w_2)$ for some $0\leq t\leq 1$. Since $R,R_1,R_2$ are tube domains, we may further assume that $$
	\im(s_1)=\im(s_2)=\im(s), \quad \im(w_1)=\im(w_2)=\im(w).
	$$
	Now consider the function $f(z)=f_{s_1,w_1,s_2,w_2}(z):=g(z s_1+(1-z)s_2, zw_1+(1-z)w_2).$ We write $$
	\begin{aligned}
		s_1=a_1+i a_2,\ \  s_2=b_1+ib_2,\ \  w_1=c_1+ic_2,\ \  w_2=d_1+i d_2,\ \  z=z_1+iz_2,
	\end{aligned}
	$$ where $a_2=b_2=\im(s)$ and $c_2=d_2=\im(w)$. Then we have
	$$
	\begin{aligned}
		z s_1+(1-z) s_2&=z_1a_1+(1-z_1)b_1+i(a_1z_2-z_2b_1+b_2),\\
		z w_1+(1-z)w_2&=z_1c_1+(1-z_1)d_1+i(c_1z_2-z_2d_1+d_2),
	\end{aligned}
	$$
	and because $(s_1,w_1)\in R_1$, $(s_2,w_2)\in R_2$, $R=\ch(R_1,R_2)$ and $R_1,R_2,R$ are open tube domains, it follows that $f(z)$ is a well-defined function that is holomorphic at least in a neighbourhood of the vertical strip $0\leq \re(z)\leq 1$.
	Now by repeatedly using the triangle inequality together with inequality $1+|ab|\leq (1+|a|)(1+|b|)$ for $a,b\in \R$, we have for $z_1=0$ the bound
	$$
	\begin{aligned}
		|f(iz_2)|&=|g(b_1+i(z_2(a_1-b_1)+b_2), d_1+i(z_2(c_1-d_1)+d_2))|\\
		&\ll (1+|b_2|)^{\alpha_2}(1+|a_1-b_1|)^{\alpha_2}(1+|d_2|)^{\beta_2}(1+|c_1-d_1|)^{\beta_2}(1+|z_2|)^{\alpha_2+\beta_2},
	\end{aligned}
	$$ and for $z_1=1$, we have
	$$
	\begin{aligned}
		|f(1+iz_2)|&=|g(a_1+i(z_2(a_1-b_1)+a_2), c_1+i(z_2(c_1-d_1)+c_2))|\\
		&\ll (1+|a_2|)^{\alpha_1}(1+|a_1-b_1|)^{\alpha_1}(1+|c_2|)^{\beta_1}(1+|c_1-d_1|)^{\beta_1}(1+|z_2|)^{\alpha_1+\beta_1}.
	\end{aligned}
	$$
	We now prove the theorem assuming that $f(z)$ is of finite order. In this case, the usual Phragm\'en-Lindel\"of principle (see \cite[Theorem 5.3]{IK}) gives for any $0\leq t\leq 1$ (remembering that $a_2=b_2$ and $c_2=d_2$):
	$$
\begin{aligned}
		|f(t+iz_2)|&\ll (1+|a_2|)^{t\alpha_1+(1-t)\alpha_2}(1+|c_2|)^{t\beta_1+(1-t)\beta_2}(1+|a_1-b_1|)^{t\alpha_1+(1-t)\alpha_2}\\
		&\cdot\quad (1+|c_1-d_1|)^{t\beta_1+(1-t)\beta_2} (1+|z_2|)^{t(\alpha_1+\beta_1)+(1-t)(\alpha_2+\beta_2)}.
\end{aligned}
	$$
	Taking $z_2=0$ and recalling that $a_2=\im(s)$, $c_2=\im(w)$ implies the bound \eqref{eq:bound in vertical strips final}. The second assertion follows from the fact that a bounded tube subdomain of $R$ can be written as the convex hull of bounded tube subdomains of $R_1,R_2$, so the quantities $|\re(s_1)-\re(s_2)|, |\re(w_1)-\re(w_2)|$ are bounded in this case. 
	
	It remains to show that $f(z)$ is of finite order.
	Let $T_1\subset R_1$, $T_2\subset R_2$ be bounded tube domains such that $(s_1,w_1)\in T_1$, $(s_2,w_2)\in T_2$ and $T_1\cap T_2\neq \emptyset$. Then Lemma \ref{le:first order in convex hulls} implies that for any $(s,w)\in T=\ch(T_1\cup T_2)$ and $\epsilon>0$, $$|g(s,w)|\ll e^{\epsilon(|\im(s)|+|\im(w)|)}.$$
	Therefore
	$$
	\begin{aligned}
		|f(z)|&\ll e^{\epsilon(|\im(z s_1+(1-z)s_2)|+|\im(z w_1+(1-z)w_2)|)}\\
		&\ll e^{\epsilon (|b_2|+|d_2|)}\cdot e^{\epsilon (|a_1|+|b_1|+|c_1|+|d_1|)|\im z|},
	\end{aligned}
	$$ which implies that $f(z)$ is of finite order.
\end{proof}

We now prove a similar result for meromorphic functions.
\begin{lemma}\label{le:bounded in vertical strips}
	Let $R_1,R_2\subset \C^2$ be open tube domains such that $R_1\cap R_2\neq\emptyset$ and $g(s,w)$ be a function that is meromorphic on $R_1,R_2$ and bounded in vertical strips there. Then it is also bounded in vertical strips on $R:=\ch(R_1,R_2).$
\end{lemma}
\begin{proof}
	Let $T\subset R$ be a compact tube subdomain of $R$. Then $T$ is bounded, so there are bounded open tube domains $T_1,T_2$ satisfying $\overline{T_1}\subset R_1$, $\overline{T_2}\subset R_2$, $T_1\cap T_2\neq\emptyset$ and $T\subset \ch(T_1,T_2)$. By Definition \ref{def:bounded meromorphic}, there are linear functions $\l_1(s,w),\dots, \l_n(s,w)$ such that $g\cdot \l_1\dots \l_n$ is holomorphic and bounded in vertical strips on $T_1\cup T_2$. The conclusion now follows from Bochner's Tube Theorem and Lemma \ref{le:bounded in vertical strips holom} applied on the holomorphic function $g\cdot \l_1\dots \l_n$.
\end{proof}

\section{Basic properties of $A(s,w)$}\label{sec:basic}

In this section, we prove the following Proposition, stating the basic properties of $A(s,w)$.
\begin{proposition}\label{prop:basic properties of A}
	The double Dirichlet series $A(s,w)$ has the following properties:
	\begin{enumerate}
		\item It satisfies the functional equation 
		\begin{equation}\label{eq:fe for A}
			A(s,w)=X_+(s)A(1-s,s+w-1/2),\quad X_+(s)=\bfrac{\pi}{8}^{s-1/2}\frac{\Gamma\bfrac{1-s}{2}}{\Gamma\bfrac{s}{2}}.
		\end{equation}
	\item It has a meromorphic continuation to the region
	$$
		R=\{\re(w)>0,\ \re(s+w)>3/2\}.
	$$
	The only pole of $\zeta(2w)A(s,w)$ in this region is a simple pole at $w=1$ with residue
	$$
		\res{w=1}A(s,w)=\frac2{3\zeta(2)}\sum_{\substack{n\geq 1\\n\odd}}\frac{\prod_{p|n}\frac{p}{p+1}}{n^{2s}}=\frac{2}{3\zeta(2)}\zeta_{(2)}(2s)P(s)=:R_1(s),
	$$
where $$
		P(s)=\prod_{p>2}\lz1-\frac{1}{p^{2s}(p+1)}\pz
$$ is an Euler product that converges absolutely for $\re(s)>0$.
	\item The function $(w-1)\zeta(2w)A(s,w)$ is polynomially bounded in vertical strips in $R$.
	\end{enumerate} 
\end{proposition}
\begin{proof}
	Recall the definition
	\begin{equation}\label{eq:def of A in Prop}
		A(s,w)=\sumstar_{\substack{d\geq 1\\d\odd}}\frac{L(s,\chi_{8d})}{d^w}.
	\end{equation}

	(1): This follows directly after applying the functional equation of $L(s,\chi_{8d})$ in \eqref{eq:def of A in Prop}, and the fact that all the occurring characters $\chi_{8d}$ are primitive modulo $8d$. We also remark that for any fixed $s$, by the polynomial boundedness of $L(s,\chi)$ in vertical strips, both sums in \eqref{eq:fe for A} are absolutely convergent if $\re(w)$ is large enough. In particular, by Lemma \ref{le:lindelof on average} (a) and partial summation, \eqref{eq:def of A in Prop} is defined for $\re(w)>\max\{1,3/2-\re(s)\}$.
	\smallskip
	
	(2): As remarked above, the definition and Lemma \ref{le:lindelof on average}, $A(s,w)$ converges absolutely (and thus is holomorphic) in the region 
	$$
		R_1=\{\re(w)>1,\ \re(s+w)>3/2\}.
	$$
We also have
\begin{equation}\label{eq:second expr for A}
	\begin{aligned}
		A(s,w)&=\sumstar_{\substack{d\geq 1\\ d\odd}}\sum_{n\geq 1}\frac{\chi_{8d}(n)}{d^wn^s}=\sum_{n\geq 1}\frac{\leg{8}{n}}{n^s}\sum_{d\geq 1}\frac{\leg{d}{4n}\mu^2(d)}{d^w}\\
		&=\sum_{n\geq 1}\frac{\leg{8}{n}L\lz w,\leg{\cdot}{4n}\pz}{n^sL\lz 2w,\leg{\cdot}{4n}^2\pz}=\frac{1}{\zeta(2w)}\sum_{n\geq 1}\frac{\leg{8}{n}L\lz w,\leg{\cdot}{4n}\pz a_{2w}(4n)}{n^s},
	\end{aligned}
\end{equation}
  where we used the notation \eqref{eq:def of a_t(n)}. By Lemma \ref{le:lindelof on average} and \eqref{eq:bound for a_t(n)} the last series converges absolutely in the region
$$
	R_2=\{\re(w)>0,\ \re(s)>1,\ \re(s+w)>3/2\},
$$ up to a pole at $w=1$ coming from $L\lz w,\leg{\cdot}{4n}\pz$ when $n$ is a square.
By Bochner's Tube Theorem, we have a meromorphic continuation of $A(s,w)$ to the region
$$
	R=\mathrm{conv}(R_1, R_2)=\{\re(w)>0,\ \re(s+w)>3/2\}.
$$
$A(s,w)$ has no pole in $R_1$ and from \eqref{eq:second expr for A}, we see that the only pole of $\zeta(2w)A(s,w)$ in $R_2$, and hence also in $R$, is at $w=1$ and comes from the terms with $n=\square$. We have
$$
\begin{aligned}
		\res{w=1}L\lz w,\leg{\cdot}{4n^2}\pz&=\res{w=1}\lz\zeta(w)\prod_{p|2n}\lz1-\frac{1}{p^w}\pz\pz=\prod_{p|2n}\lz1-\frac{1}{p}\pz\\
		&=a_1(2n)^{-1},
\end{aligned}
$$
so we obtain
$$
	\res{w=1}A(s,w)=\frac{1}{\zeta(2)}\sum_{\substack{n\geq 1\\n\odd}}\frac{a_2(2n)/a_1(2n)}{n^{2s}}=\frac2{3\zeta(2)}\sum_{\substack{n\geq 1\\n\odd}}\frac{\prod_{p|n}\frac{p}{p+1}}{n^{2s}}.
$$
Writing the last Dirichlet series as an Euler product and using $\lz1+\frac{1}{p^{2s}-1}\pz^{-1}=\lz1-\frac{1}{p^{2s}}\pz$, we obtain
$$
\begin{aligned}
		\sum_{\substack{n\geq 1\\n\odd}}\frac{\prod_{p|n}\frac{p}{p+1}}{n^{2s}}&=\prod_{p>2}\lz 1+\frac{p}{p+1}\cdot\frac{1}{p^{2s}-1}\pz=\prod_{p>2}\lz 1+\frac{1}{p^{2s}-1}-\frac{1}{(p+1)(p^{2s}-1)}\pz\\
		&=\prod_{p>2}\lz1+\frac{1}{p^{2s}-1}\pz\prod_{p>2}\lz1-\frac{1}{(p+1)(p^{2s}-1)}+\frac{1}{(p+1)p^{2s}(p^{2s}-1)}\pz\\
		&=\zeta_{(2)}(2s)P(s),
\end{aligned}
$$ where $$
P(s)=\prod_{p>2}\lz1-\frac{1}{(p+1)(p^{2s}-1)}+\frac{1}{(p+1)p^{2s}(p^{2s}-1)}\pz=\prod_{p>2}\lz1-\frac{1}{p^{2s}(p+1)}\pz
$$ is an Euler product which is absolutely convergent for $\re(s)>0$.

(3): This follows from Lemma \ref{le:bounded in vertical strips} and the fact that $A(s,w)$ is polynomially bounded in vertical strips in both $R_1$ and $R_2$.
\end{proof}

\section{Meromorphic continuation of $A(s,w)$ and the functional equation in $w$}\label{sec:fe in w}
In this section, we provide further meromorphic continuation of $A(s,w)$ using the ``functional equation in $w$". This is slightly more complicated than the ``functional equation in $s$" from Proposition \ref{prop:basic properties of A} (1) due to the presence of non-primitive characters. In this case, we have a functional equation provided by \eqref{eq:general fe}, which gives a relation between $A(s,w)$ and a different double Dirichlet series, whose coefficients are Gauss sums.

\begin{lemma}\label{le:fe in w} We have
	\begin{equation}\label{eq:fe in w}
		\begin{aligned}
			A(s,w)&=\frac{Y_+(w)}{\zeta(2w)}B_+(s+w-\tfrac12,1-w;2w)+\frac{Y_-(w)}{\zeta(2w)}B_-(s+w-\tfrac12,1-w;2w),
		\end{aligned}
	\end{equation} where
	\begin{equation}\label{eq:B1}
		\begin{aligned}
			B_{\pm}(s,w;t)&=\sum_{n\equiv\pm 1\mod 4}\frac{K\lz w,\leg{\cdot}{4n}\pz a_t(4n)\psi_2(n)}{n^s}\\
			&=\frac12\sum_{\substack{n\equiv\pm 1\mod 4\\k\geq 1}}\frac{\tau\lz\leg{\cdot}{4n},k\pz a_t(4n)\psi_2(n)}{n^{s+1/2}k^w}.
		\end{aligned}
	\end{equation}
\end{lemma}

\begin{proof} The proof follows directly from the last expression in \eqref{eq:second expr for A} and Proposition \ref{prop:general fe}. 
\end{proof}
In order to obtain a meromorphic continuation of $A(s,w)$, our goal now is to determine the regions of absolute convergence or meromorphic continuation of $B_{\pm}(s,w;t)$.

One such region can be obtained from the first expression in \eqref{eq:B1} and Lemma \ref{le:lindelof for K}. A second region will be obtained after expanding the Dirichlet series and exchanging summations. For this, it will be useful to write $B_\pm(s,w;t)$ in terms of Dirichlet series whose coefficients are multiplicative, which is achieved in the following lemma by replacing the Gauss sums with $G\lz\leg{\cdot}{n},d\pz$ and capturing the congruence conditions using Dirichlet characters.
\begin{lemma}\label{le:rewriting B}
	We have
\begin{equation}\label{eq:express B}
\begin{aligned}
			B_{\pm}(s,w;t)&=\frac{\epsilon_\pm a_t(2)}{2\cdot 4^w}\lz B\lz s,w;t,\psi_2,\psi_0\pz\pm B\lz s,w;t,\psi_{-2},\psi_0\pz\pz\\
	&-\frac{\epsilon_\pm a_t(2)}{2\cdot 2^w}\lz B(s,w;t,\psi_1,\psi_1)\pm B(s,w;t,\psi_{-1},\psi_1) \pz,
\end{aligned}
\end{equation} where $\epsilon_+=1$, $\epsilon_-=i$, and
$$
	B(s,w;t,\psi,\psi')=\sum_{n,k\geq 1}\frac{G\lz\leg{\cdot}{n},k\pz a_t(n)\psi(n)\psi'(k)}{n^{s+1/2}k^w}.
$$
\end{lemma}
\begin{proof}
 By \eqref{eq:Gauss sums with 4} and \eqref{eq:gauss sums second coord} (and noting that the sum runs only over odd $n$), we have
$$
\begin{aligned}
		B_{\pm}(s,w;t)&=\sum_{n\equiv\pm 1\mod 4}\sum_{\substack{k\geq 1\\k\equiv 0\mod 4}}\frac{\tau\lz\leg{\cdot}{n},k\pz a_{t}(4n)\psi_2(n)}{n^{s+1/2}k^w}\\
		&-\sum_{n\equiv\pm1\mod 4}\sum_{\substack{k\geq 1\\k\equiv 2\mod 4}}\frac{\tau\lz\leg{\cdot}{n},k\pz a_{t}(4n)\psi_2(n)}{n^{s+1/2}k^w}\\
		&=\frac{a_t(2)}{4^w}\sum_{\substack{n\equiv\pm 1\mod 4\\k\geq 1}}\frac{\tau\lz\leg{\cdot}{n},k\pz a_t(n)\psi_2(n)}{n^{s+1/2}k^w}\\
		&-\frac{a_t(2)}{2^w}\sum_{\substack{n\equiv\pm 1\mod 4\\k\geq 1}}\frac{\tau\lz\leg{\cdot}{n},k\pz a_t(n)\psi_1(k)}{n^{s+1/2}k^w}.
\end{aligned}
$$
Now letting $\epsilon_+=1$ and $\epsilon_-=i$, \eqref{eq:tau as G} gives
$$
\begin{aligned}
		B_{\pm}(s,w;t)&=\frac{\epsilon_\pm a_t(2)}{4^w}\sum_{\substack{n\equiv\pm 1\mod 4\\ k\geq 1}}\frac{G\lz\leg{\cdot}{n},k\pz a_t(n)\psi_2(n)}{n^{s+1/2}k^w}\\
		&-\frac{\epsilon_\pm a_t(2)}{2^w}\sum_{\substack{n\equiv \pm1\mod 4\\k\geq 1}}\frac{G\lz\leg{\cdot}{n},k\pz a_t(n)\psi_1(k)}{n^{s+1/2}k^w}\\
		&=\frac{\epsilon_\pm a_t(2)}{4^w}\sum_{n,k\geq 1}\frac{G\lz\leg{\cdot}{n},k\pz a_t(n)\psi_2(n)}{n^{s+1/2}k^w}\lz\frac{\psi_1(n)\pm\psi_{-1}(n)}{2}\pz\\	
		&-\frac{\epsilon_\pm a_t(2)}{2^w}\sum_{n,k\geq 1}\frac{G\lz\leg{\cdot}{n},k\pz a_t(n)\psi_1(k)}{n^{s+1/2}k^w}\lz\frac{\psi_1(n)\pm\psi_{-1}(n)}{2}\pz\\
		&=\frac{\epsilon_\pm a_t(2)}{2\cdot 4^w}\lz B\lz s,w;t,\psi_2,\psi_0\pz\pm B\lz s,w;t,\psi_{-2},\psi_0\pz\pz\\
		&-\frac{\epsilon_\pm a_t(2)}{2\cdot 2^w}\lz B(s,w;t,\psi_1,\psi_1)\pm B(s,w;t,\psi_{-1},\psi_1) \pz
\end{aligned}
$$
\end{proof}

Now we write 
\begin{equation}\label{eq:definition of B(s,w;t,psi,psi')}
	B(s,w;t,\psi,\psi')=\sum_{k\geq 1}\frac{\psi'(k)}{k^w}\sum_{n\geq 1}\frac{G\lz\leg{\cdot}{n},k\pz a_t(n)\psi(n)}{n^{s+1/2}},
\end{equation} and notice that the coefficients of the inner sum are multiplicative in $n$. Before expanding it into an Euler product, we prove the following auxiliary lemma.

\begin{lemma}\label{le:euler product}
	For any $N\in \N$, $s,t\in\C$ with $\re(s),\re(t)>1$, and Dirichlet character $\chi$, we have
	\begin{equation}\label{eq:euler product}
			\prod_p\lz1+\frac{\chi(p)}{\lz p^s+\chi(p)\pz\lz p^t-1\pz}\pz=\prod_{j=1}^N\frac{L(s+jt,\chi)}{L(2s+2jt,\chi^2)}\cdot E_N(s,t;\chi),
	\end{equation}
where $E_N(s,t;\chi)$ is an Euler product that satisifes the recurrence relation
\begin{equation}\label{eq:def of E}
\begin{aligned}
		E_1(s,t;\chi)&=\prod_p\lz1+\frac{\chi(p)p^s-\chi(p)^2p^t+\chi(p)^2}{\lz p^{s+t}+\chi(p)\pz\lz p^s+\chi(p)\pz\lz p^t-1\pz}\pz,\\
		E_{N+1}(s,t;\chi)&=E_N(s,t;\chi)\prod_p\lz1+\frac{\chi(p)}{p^{s+(N+1)t}}\pz^{-1},\qquad\hbox{$N\geq 1$.}
\end{aligned}
\end{equation}
Moreover, for every $N\geq 1$, we have
$$
	E_N(s,t;\chi)=\prod_p\lz1+\frac{\chi(p)p^{Ns+\frac{N(N-1)}{2}t}}{(p^t-1)\prod\limits_{0\leq j\leq N}\lz p^{s+jt}+\chi(p)\pz}+O_N\lz\frac{1}{p^{2s+t}}\pz\pz.
$$
\end{lemma}
\begin{proof}
	We proceed by induction over $N$.
	
	For $N=1$, we write the LHS of \eqref{eq:euler product} as
	$$
\begin{aligned}
			&\prod_p\lz1+\frac{\chi(p)}{p^{s+t}}+\frac{\chi(p)}{(p^s+\chi(p))(p^t-1)}-\frac{\chi(p)}{p^{s+t}}\pz\\
			&=\prod_p\lz1+\frac{\chi(p)}{p^{s+t}}+\frac{\chi(p)p^{s+t}-\chi(p)(p^s+\chi(p))(p^t-1)}{p^{s+t}(p^s+\chi(p))(p^t-1)}\pz\\
			&=\frac{L(s+t,\chi)}{L(2s+2t,\chi^2)}\prod_p\lz1+\frac{\chi(p)p^{s+t}-\chi(p)(p^s+\chi(p))(p^t-1)}{p^{s+t}(p^s+\chi(p))(p^t-1)}\lz1+\frac{\chi(p)}{p^{s+t}}\pz^{-1}\pz\\
			&=\frac{L(s+t,\chi)}{L(2s+2t,\chi^2)}E_1(s,t;\chi),
\end{aligned}
	$$ where
$$
	E_1(s,t;\chi)=\prod_p\lz1+\frac{\chi(p)p^s-\chi(p)^2 p^t+\chi(p)^2}{\lz p^{s+t}+\chi(p)\pz(p^s+\chi(p))(p^t-1)}\pz.
$$ We also have
$$
\begin{aligned}
		E_1(s,t;\chi)&=\prod_{p}\lz1+\frac{\chi(p)p^s+O(p^t)}{\lz p^{s+t}+\chi(p)\pz\lz p^s+\chi(p)\pz\lz p^t-1\pz}\pz\\
	&=
	\prod_p\lz1+\frac{\chi(p)p^s}{(p^{s+t}+\chi(p))(p^s+\chi(p))(p^t-1)}+O\bfrac{1}{p^{2s+t}}\pz,
\end{aligned}
$$ which finishes the proof when $N=1$.

Now assume that $P(s,t;\chi)=\prod_{j=1}^N\frac{L(s+jt,\chi)}{L(2s+2jt,\chi^2)} \cdot E_N(s,t;\chi)$, where
$$
	E_N(s,t;\chi)=\prod_p\lz 1+\frac{\chi(p)p^{Ns+\frac{N(N-1)}{2}t}}{\lz p^t-1\pz\prod\limits_{0\leq j\leq N}\lz p^{s+jt}+\chi(p)\pz}+O\lz\frac{1}{p^{2s+t}}\pz\pz.
$$
Then also
$P(s,t;\chi)=\prod_{j=1}^{N+1}\frac{L(s+jt,\chi)}{L(2s+2jt,\chi^2)} \cdot E_{N+1}(s,t;\chi)$, where
$$
\begin{aligned}
		&E_{N+1}(s,t;\chi)=E_N(s,t;\chi)\prod_p\lz1+\frac{\chi(p)}{p^{s+(N+1)t}}\pz^{-1}\\
		&=\prod_p\lz 1+\frac{\chi(p)p^{Ns+\frac{N(N-1)}{2}t}}{\lz p^t-1\pz\prod\limits_{0\leq j\leq N}\lz p^{s+jt}+\chi(p)\pz}+O\lz\frac{1}{p^{2s+t}}\pz\pz\\
		&\qquad\cdot\qquad\lz1+\frac{\chi(p)}{p^{s+(N+1)t}}\pz^{-1}\\
		&=\prod_p\Bigg\{1+O\bfrac{1}{p^{2s+t}}+\lz\frac{\chi(p)p^{Ns+\frac{N(N-1)}{2}t}}{\lz p^t-1\pz\prod\limits_{0\leq j\leq N}\lz p^{s+jt}+\chi(p)\pz}-\frac{\chi(p)}{p^{s+(N+1)t}}\pz\\
		&\qquad\qquad\qquad\qquad\qquad\cdot\qquad\lz1+\frac{\chi(p)}{p^{s+(N+1)t}}\pz^{-1}\Bigg\}.
	\end{aligned}
	$$
We have
$$
\begin{aligned}
		&\lz\frac{\chi(p)p^{Ns+\frac{N(N-1)}{2}t}}{\lz p^t-1\pz\prod\limits_{0\leq j\leq N}\lz p^{s+jt}+\chi(p)\pz}-\frac{\chi(p)}{p^{s+(N+1)t}}\pz\lz1+\frac{\chi(p)}{p^{s+(N+1)t}}\pz^{-1}\\
		&=\frac{\chi(p)p^{(N+1)s+\frac{N^2+N+2}{2}t}-\chi(p)(p^t-1)\prod\limits_{0\leq j\leq N}\lz p^{s+jt}+\chi(p)\pz}{(p^t-1)\prod\limits_{0\leq j\leq N}\lz p^{s+jt}+\chi(p)\pz\cdot p^{s+(N+1)t}}\lz1+\frac{\chi(p)}{p^{s+(N+1)t}}\pz^{-1}\\
		&=\frac{\chi(p)p^{(N+1)s+\frac{N(N+1)}{2}t}+O\lz p^{Ns+\frac{N^2+N+2}{2}t}\pz}{(p^t-1)\prod\limits_{0\leq j\leq N+1}\lz p^{s+jt}+\chi(p)\pz}\\
		&=\frac{\chi(p)p^{(N+1)s+\frac{N(N+1)}{2}t}}{(p^t-1)\prod\limits_{0\leq j\leq N+1}\lz p^{s+jt}+\chi(p)\pz}+O\bfrac{1}{p^{2s+(N+1)t}},
\end{aligned}
$$
so
$$
\begin{aligned}
		E_{N+1}(s,t;\chi)=\prod_p\lz
		1+\frac{\chi(p)p^{(N+1)s+\frac{N(N+1)}{2}t}}{(p^t-1)\prod\limits_{0\leq j\leq N+1}\lz p^{s+jt}+\chi(p)\pz}+O\bfrac{1}{p^{2s+t}}
		\pz
\end{aligned}
$$
\end{proof}

\begin{lemma}\label{le:continuation}
For any $N\in \N$, $\re(s),\re(t)>0$, and character $\psi$ with $\psi(2)=0$, we have
	\begin{equation}\label{eq:Euler product for D}
\begin{aligned}
			\sum_{n\geq 1}\frac{G\lz\leg{\cdot}{n},k\pz\psi(n)a_t(n)}{n^{s+1/2}}&=\prod_{j=0}^N\frac{L\lz s+jt,\leg{4k}{\cdot}\psi\pz}{L\lz 2s+2jt,\leg{4k}{\cdot}^2\psi^2\pz}\cdot E_N\lz s,t;\lz\tfrac{4k}{\cdot}\pz\psi\pz
		P(s,t,k;\psi),
\end{aligned}
	\end{equation}
where $E_N\lz s,t;\leg{4k}{\cdot}\psi\pz$ is as in Lemma \ref{le:euler product}, so it satisfies
$$
	E_N\lz s,t;\leg{4k}{\cdot}\psi\pz=\prod_p\lz1+O_N\lz\frac{1}{p^{s+(N+1)t}}+\frac{1}{p^{2s+t}}\pz\pz,
$$
and 
\begin{equation}\label{eq:def of P}
	P(s,t,k;\psi)=\prod_{p|k}\lz 1+a_t(p)\sum_{j\geq 1}\frac{G\lz\leg{\cdot}{p^j},k\pz\psi(p^j)}{p^{js+j/2}}\pz
\end{equation}
 is a finite Euler product, which for any $\delta,\epsilon>0$, $\re(t)>\delta$, and $\re(s)>1/2+\delta$ satisfies $|P(s,t,k;\psi)|\ll_{\epsilon,\delta} k^\epsilon$.
\end{lemma}

\begin{remark}\label{rem:square-freeness}
	The Euler product in \eqref{eq:Euler product for D} is somewhat complicated, mainly due to the presence of the factors $a_t(n)$ on the left-hand side. These factors appeared in \eqref{eq:second expr for A} because of $L\lz 2w,\lz\tfrac{\cdot}{4n}\pz^2\pz$ in the denominator, which in turn came from the square-free condition on $d$. If we considered the family where $d$ runs over all positive integers, the factors $a_t(n)$ would not be present and the product $P_{p\nmid k}$ in \eqref{eq:P not dividing k} below would simply be $\frac{L\lz s,\leg{4k}{\cdot}\psi\pz}{L\lz 2s,\leg{4k}{\cdot}^2\psi^2\pz}$. In particular, there would be no extra pole from the factor $\frac{L\lz s+t,\leg{4k}{\cdot}\psi\pz}{L\lz 2s+2t,\leg{4k}{\cdot}^2\psi^2\pz}$, which is where the lower order terms come from.
	 
	We would then recover Theorem \ref{thm:moment} without the lower order terms, which corresponds to the result of Gao and Zhao \cite{gazh_first_moment}.
\end{remark}
\begin{proof}
	Since all the coefficients are multiplicative in $n$, we have
$$
\begin{aligned}
			\sum_{n\geq 1}\frac{G\lz\leg{\cdot}{n},k\pz\psi(n)a_t(n)}{n^{s+1/2}}&=\prod_p\lz1+a_t(p)\sum_{j\geq 1}\frac{G\lz\leg{\cdot}{p^j},k\pz\psi(p^j)}{p^{js+j/2}}\pz\\
			&=:P_{p|k}P_{p\nmid k},
\end{aligned}
$$
where $P_{p|k},P_{p\nmid k}$ correspond to the factors with $p|k$, $p\nmid k$, respectively. Since $\psi(2)=0$, we may assume from now on that $p$ denotes odd primes.

For $p\nmid k$, \eqref{eq:Gauss sum evaluation} implies that $G\lz \leg{\cdot}{p},k\pz=\leg{k}{p}p^{1/2}$, and $G\lz \leg{\cdot}{p^j},k\pz=0$ for $j\geq 2$. We have $\leg{k}{p}=\leg{4k}{p}$ and $\leg{4k}{\cdot}$ is a well-defined Dirichlet character modulo $4k$, so
\begin{equation}\label{eq:P not dividing k}
\begin{aligned}
		P_{p\nmid k}&=\prod_{p\nmid k}\lz 1+\lz1-\frac{1}{p^t}\pz^{-1}\frac{\leg{4k}{p}\psi(p)}{p^s}\pz=\prod_p\lz1+\frac{\leg{4k}{p}\psi(p)}{p^s}+\frac{\leg{4k}{p}\psi(p)}{p^s(p^t-1)}\pz\\
		&=\prod_p\lz1+\frac{\leg{4k}{p}\psi(p)}{p^s}\pz\prod_p\lz1+\frac{\leg{4k}{p}\psi(p)}{\lz p^s+\leg{4k}{p}\psi(p)\pz\lz p^t-1\pz}\pz\\
		&=\frac{L\lz s,\leg{4k}{\cdot}\psi\pz}{L\lz 2s,\leg{4k}{\cdot}^2\psi^2\pz}\prod_p\lz1+\frac{\leg{4k}{p}\psi(p)}{\lz p^s+\leg{4k}{p}\psi(p)\pz\lz p^t-1\pz}\pz.
\end{aligned}
\end{equation} Now Lemma \ref{le:euler product} gives that for any positive integer $N$, we have
$$
		P_{p\nmid k}=\prod_{j=0}^N\frac{L\lz s+jt,\leg{4k}{\cdot}\psi\pz}{L\lz 2s+2jt,\leg{4k}{\cdot}^2\psi^2\pz} \cdot E_N\lz s,t; \leg{4k}{\cdot}\psi\pz,
$$
where $E_N\lz s,t; \leg{4k}{\cdot}\psi\pz$ is an Euler product that satisfies
$$
	E_N\lz s,t; \leg{4k}{\cdot}\psi\pz=\prod_p\lz1+\frac{\leg{4k}{p}\psi(p)p^{Ns+\frac{N(N-1)}{2}t}}{(p^t-1)\prod\limits_{0\leq j\leq N}\lz p^{s+jt}+\leg{4k}{p}\psi(p)\pz}+O\lz\frac{1}{p^{2s+t}}\pz\pz,
$$ and we note that the second term in the last parenthesis is $\ll_N p^{-s-(N+1)t}$. This proves part of the Lemma after setting $P(s,t,k;\psi)=P_{p|k}$.

It remains to bound $P(s,t,k;\psi)$. Observe that this is a finite product and each factor is a finite sum by \eqref{eq:Gauss sum evaluation}. Now the trivial bound gives $\lab G\lz\leg{\cdot}{p^j},k\pz\rab\leq p^{j}.$ For any fixed $\delta>0$ and $\re(t)>\delta$, we also have $a_t(p)\ll_\delta 1.$ Therefore for $\re(s)>1/2+\delta$, and a large enough constant $c$ (depending on $\delta$), we have for any $\epsilon>0$
$$
	\begin{aligned}
		|P(s,t,k;\psi)|&=\lab\prod_{p|k}\lz1+a_{t}(p)\sum_{j\geq 1}\frac{G\lz\leg{\cdot}{p^j},k\pz}{p^{js+j/2}}\pz\rab\\
		&\ll\lab\prod_{p|k}\lz1+|a_t(p)|\sum_{j\geq 1}\frac{1}{p^{j(\re(s)-1/2)}}\pz\rab\leq c^{\omega(k)}\ll_{\epsilon,\delta} k^\epsilon.
	\end{aligned}
$$

\end{proof}

\begin{proposition}\label{prop:continuation of B}
		Assume the truth of the generalized Lindel\"of Hypothesis. Then the functions $B_\pm(s,w;t)$ satisfy the following.
		\begin{enumerate}
			\item They have a meromorphic continuation to the region
$$
				S=\{\re(s)>1/2,\ \re(s+w)>3/2,\ \re(t)>0\}.
$$
			\item Their only potential poles in $S$ are at $w=1$ and $s+tj=1$ for $j\in \Z_{\geq 0}$.
			
			\item  They are polynomially bouded in vertical strips in $S$.
		\end{enumerate} 
\end{proposition}

\begin{proof}
(1): From the definition, we have
	\begin{equation}\label{eq:B in S1}
		B_\pm(s,w;t)=\sum_{n\equiv\pm1\mod 4}\frac{K\lz w,\leg{\cdot}{4n}\pz a_t(4n)}{n^s}.
	\end{equation}
By Lemma \ref{le:lindelof for K}, \eqref{eq:bound for a_t(n)} and partial summation, this defines a holomorphic function in the region
$$
	S_1=\{\re(s)>1,\ \re(s+w)>3/2,\ \re(t)>0\}.
$$
To obtain a meromorphic continuation, Lemma \ref{le:rewriting B} implies that it is enough to find a meromorphic continuation of each $$
	B(s,w;t,\psi,\psi')=\sum_{k\geq 1}\frac{\psi'(k)}{k^w}\sum_{n\geq 1}\frac{G\lz\leg{\cdot}{n},k\pz a_t(n)\psi(n)}{n^{s+1/2}},
$$ where $\psi$ is a character that satisfies $\psi(2)=0$. For any $\delta>0$ and $N\in\N$, Lemma \ref{le:continuation} and the Lindel\"of Hypothesis implies that $B(s,w;t,\psi,\psi')$, and hence also $B_{\pm}(s,w;t)$, has a meromorphic continuation to the region
$$
	S_2(N,\delta)=\{\re(w)>1,\ \re(s)>1/2+\delta,\ \re(t)>\delta,\ \re(s+(N+1)t)>1\}.
$$
Since for any $\delta>0$, the intersection $\bigcap_{N=1}^\infty S_2(N,\delta)$ is non-empty, we can obtain a meromorphic continuation to the union of these regions, which is given by
\begin{equation}\label{eq:region N}
	S_2=\bigcup\limits_{\delta>0,N\in \N}S_2(N,\delta)=\{\re(w)>1,\ \re(s)>1/2,\  \re(t)>0\}.
\end{equation} 
Bochner's Tube Theorem provides a meromorphic continuation to 
$$
	S=\mathrm{conv}(S_1\cup S_2)=\{\re(s)>1/2,\ \re(s+w)>3/2,\ \re(t)>0\}.
$$

\smallskip

(2): Now we determine the possible poles of $B_\pm(s,w;t)$. In $S_1$, we have \eqref{eq:B in S1}, so the only poles in this region are at the poles of $K\lz w,\leg{\cdot}{4n}\pz$, which can occur at $w=1$ if $n$ is a square.

For any $N \in \N$, $\delta>0$, for $(s,w,t)\in S_2(N,\delta)$, we have
\begin{equation}\label{eq:B in S2}
\begin{aligned}
		&B(s,w;t,\psi,\psi')\\
		&=\sum_{k\geq 1}\frac{\psi'(k)}{k^w}\prod_{j=0}^N\frac{L(s+tj,\leg{4k}{\cdot}\psi)}{L\lz 2s+2tj,\leg{4k}{\cdot}^2\psi^2\pz}\cdot E_N\lz s,t;\leg{4k}{\cdot}\psi\pz P(s,t,k;\psi),
\end{aligned}
\end{equation} and Lemma \ref{le:rewriting B} implies that the potential poles of $B_{\pm}(s,w;t)$ come only from the potential poles of $B(s,w;t,\psi,\psi')$. Now in \eqref{eq:B in S2} the poles arise only from the poles of the summands, as otherwise the sum over $k$ is absolutely convergent in $S_2(N,\delta).$ The products $E_N$ and $P$ are holomorphic throughout the whole region $S_2(N,\delta)$, so the poles arise only from the possible poles of $L\lz s+tj,\leg{4k}{\cdot}\psi\pz$ and zeros of $L\lz 2s+2tj,\leg{4k}{\cdot}^2\psi^2\pz$. The denominators have zeros at the points $2s+2tj=\rho$, where $0\leq j\leq N$, and $\rho$ is a zero of $\zeta_{(2k)}(s)$. These satisfy $\re(s+tj)<1/2$, so they are outside $S_2(N,\delta).$

The remaining poles come from the terms where $\leg{4k}{\cdot}\psi$ is a trivial character, which can occur only if $\psi=\psi_1$ and $k=\square$, or if $\psi=\psi_2$ and $k=2\cdot\square$, and the poles are at $s+tj=1$.

\smallskip

(3):
Let $U\subset T$ be a compact tube subdomain of $T$, and let $\delta>0$ be small and $N$ be large enough so that $U\subset \ch\lz S_1,\bigcup\limits_{j=1}^N S(j,\delta)\pz$, we show polynomial boundedness on the larger set. We see from \eqref{eq:B in S1} that $B_{\pm}(s,w;t)$ are polynomially bounded in vertical strips in $S_1$, and from \eqref{eq:B in S2} and Lemma \ref{le:rewriting B} that the same holds in each $S(j,\delta)$, so also in $\bigcup\limits_{j=1}^N S(j,\delta)$. The conclusion now follows from Lemma \ref{le:bounded in vertical strips}.
\end{proof}

We now obtain the meromorphic continuation and location of poles of $A(s,w)$, which gives the first part of Theorem \ref{thm:MDS}.

\begin{proposition}\label{prop:continuation of A final}
	The double Dirichlet series $A(s,w)$ has a meromorphic continuation to the region
	$$
		T=\{\re(s+w)>1/2,\ \re(s+2w)>1,\ \re(w)>0\}.
	$$
The possible poles of $\zeta(2w)\zeta(2s+2w-1)A(s,w)$ in this region occur at $w=1$, $s+w=3/2$, $s+(2j+1)w=3/2$ and $2js+(2j+1)w=j+1$ for $j\in\N$. It is polynomially bounded in vertical strips in $T$.
\end{proposition}
\begin{proof}
	By Proposition \ref{prop:basic properties of A}, $A(s,w)$ is defined and meromorphic throughout $R_1:=\{\re(w)>0,\ \re(s+w)>3/2\}$, and the only pole of $\zeta(2w)A(s,w)$ in $R$ is at $w=1$. By Lemma \ref{le:fe in w} and Proposition \ref{prop:continuation of B}, we obtain a meromorphic continuation of $A(s,w)$ to the region $R_2 := \{\re(s)>1,\re(w)>0\}$, and the potential poles of $\zeta(2w)A(s,w)$ in this region are at $s+w=3/2$, $s+(2j+1)w=3/2$ for $j\in\N$ (by the remarks under Proposition \ref{prop:general fe}, the poles of $Y\pm(s)$ do not give rise to new poles of $A(s,w)$).
	
	By Bochner's Tube Theorem, $\zeta(2w)A(s,w)$ has a meromorphic continuation to
	$$
		T_1:=\ch{(R_1,R_2)}=\{\re(s+w)>1,\ \re(w)>0\},
	$$ with only possible poles at $w=1$, and $s+(2j+1)w=3/2$ for $j\in\Z_{\geq 0}$. The same holds for $\zeta(2s+2w-1)\zeta(2w)A(s,w)$, since the extra pole of $\zeta(2s+2w-1)$ lies outside of $T_1$.

	Using the functional equation in Proposition \ref{prop:basic properties of A}, denoting $\sigma$ the affine transform $\sigma:(s,w)\mapsto(1-s,s+w-1/2)$ and noting that $\sigma(2w)=2s+2w-1$, we obtain a meromorphic continuation of $\zeta(2w)\zeta(2s+2w-1)A(s,w)$ to the region 
	$$
		T_2=\sigma(T_1)=\{\re(w)>1/2,\ \re(s+w)>1/2\}
	$$ with potential poles in $T_2$ at $s+w=3/2$, and $2js+(2j+1)w=j+1$ for $j\in\Z_{\geq 0}$.

The final result follows from a last application of Bochner's Tube Theorem, which provides meromorphic continuation to 
$$
	T=\ch(T_1,T_2)=\{\re(s+w)>1/2,\ \re(s+2w)>1,\ \re(w)>0\},
$$ and the possible poles of $\zeta(2w)\zeta(2s+2w-1)A(s,w)$ are at $w=1$, $s+w=3/2$, $s+(2j+1)w=3/2$ and $2js+(2j+1)w=j+1$ for $j=1,2,\dots$. Polynomial boundedness follows from Lemma \ref{le:bounded in vertical strips} and the polynomial boundedness of $A(s,w)$ and $B_\pm(s,w;2w)$ in the corresponding regions, as established in propositions \ref{prop:basic properties of A} and \ref{prop:continuation of B}. 
\end{proof}

\begin{remark}\label{rem:no lindelof}
	If we did not assume the generalized Lindel\"of hypothesis, we would be able to take the union with $N=3$ in \eqref{eq:region N}, since GLH holds on average for the product of 4 real Dirichlet $L$-functions by Heath-Brown's quadratic large sieve \cite{HB_quadratic_sieve}. We would then obtain Proposition \ref{prop:continuation of A final} with the region $$T = \{\re(s+w)>\tfrac12, \re(w)>0\, \re(8s+9w)>5, \re(s+9w)>\tfrac32,\re(s+2w)>\tfrac{10}{9}\}.$$ This would lead to Theorem \ref{thm:moment} with a slightly weaker error. For instance at the central point $s=1/2$, we would have the extra condition $\re(w)>11/36$, giving the error $O(X^{11/36+\epsilon})$. However, this is still more than enough to unconditionally recover the error term $O(X^{1/2+\epsilon})$, and even the lower-order terms under RH.
\end{remark}

To finish the proof of Theorem \ref{thm:MDS}, we need to compute the claimed residues, which we do in the following section.

\section{computing the residues}\label{sec:residues}
To finish the proof of Theorem \ref{thm:MDS}, we compute the residues of $A(s,w)$. 

\begin{proposition}\label{prop:residues}
	Assume that $0<\re(s)<3/2$, $s\neq 1/2$. Then
	\begin{enumerate}
		\item 
			We have $\res{w=1}A(s,w)=\frac{2}{3\zeta(2)}\zeta_{(2)}(2s)P(s)$,
			where
			$$
			P(s)=\prod_{p>2}\lz1-\frac{1}{p^{2s}(p+1)}\pz.
			$$
		\item We have
		$\res{w=3/2-s}A(s,w)=X_+(s)\res{w=1}A(1-s,w).$
		\item We have
$\res{w=1/2-s/3} A(s,w)=Y(s)\zeta(2s)Q_2(s)Q(s)$, where
$$
	\begin{aligned}
		Y(s)&:=Y_+\bfrac{2s}{3}\lz Y_+\bfrac{3-2s}{6}+iY_-\bfrac{3-2s}{6}\pz,\\
		Q_2(s)&:=\lz1-\frac{1}{2^{2s}}\pz\frac{2^{\frac{2s}3}(4^{\frac12+\frac s3}-2)}{3(2^{\frac32+s}-2^{\frac12+\frac s3})(2^{\frac{4s}3}-2^{\frac{2s}3+1}-2)},\\
		Q(s)&:=\frac{E\bfrac{2s}{3}}{\zeta\bfrac{4s}{3}\zeta(2)}\prod_{p>2}\lz1+\frac{p^{1+4s/3}+p^{1+2s/3}-p-p^{2s}}{p^{2s/3}(p^{1+2s/3}+p-p^{4s/3})(p^{1+4s/3}-1)}\pz,
	\end{aligned}
$$ and 
$$
	E(s)=\prod_p\lz 1+\frac{p^s-p^{1-s}+1}{(p+1)(p^s+1)(p^{1-s}-1)}\pz.
$$
\item We have
		$\res{w=2/3-2s/3} A(s,w)=X_+(s)\res{w=1/2-s/3}A(1-s,w).$
		\end{enumerate}
\end{proposition}

\begin{proof}[Proof of Proposition \ref{prop:residues} (1), (2) and (4)]
	
(1) and (2): the residue at $w=1$ was computed in Proposition \ref{prop:basic properties of A} (b). The residue at $w=3/2-s$ comes from the functional equation from Proposition \ref{prop:basic properties of A} (a) as shown in the following computation, where we use the functional equation and change variables $v=s+w-1/2$:
\begin{equation}\label{eq:residue and functional equation}
\begin{aligned}
		\res{w=\frac32-s}A(s,w)&=\lim_{w\rightarrow \frac32-s}(w+s-\tfrac32)A(s,w)\\
		&=\lim_{w\rightarrow \frac32-s}(w+s-\tfrac32)X_+(s)A(1-s,s+w-\tfrac12)\\
		&=\lim_{v\rightarrow 1}(v-1)X_+(s)A(1-s,v)\\
		&=X_+(s)\res{v=1}A(1-s,v).
\end{aligned}
\end{equation}

(4): The functional equation $A(s,w)=X_+(s)A(1-s,s+w-1/2)$ transforms the pole at $w=1/2-s/3$ into the pole at $s+w-1/2=1/2-1/3+s/3$, so at $w=2/3-2s/3$, and a computation similar to that in \eqref{eq:residue and functional equation} gives
$$
	\res{w=2/3-2s/3}A(s,w)=X_+(s)\res{w=1/2-s/3} A(1-s,w).
$$ 
\end{proof} 

It remains to compute the residue at $w=1/2-s/3$, which arises from the poles of $B_\pm(s,w;t)$ at $s+t=1$ and the functional equation \eqref{eq:fe in w}.

We will use the following notation:
	$$
		\begin{aligned}
			R(s,w;\psi,\psi')&:=\res{t=1-s}B(s,w;t,\psi,\psi'),\\
			R_\pm(s,w)&:=\res{t=1-s}B_{\pm}(s,w;t),\\
			R(s)&:=\res{w=1/2-s/3}A(s,w).
		\end{aligned}
	$$

We note that $R(s,w;\psi_{-1},\psi')=R(s,w;\psi_{-2},\psi')=0$, because in this case $\leg{4k}{\cdot}\psi$ from \eqref{eq:Euler product for D} is non-principal for any $k\geq 1$, so \eqref{eq:express B} gives
\begin{equation}\label{eq:R+ as R(s,w,psi,psi')}
		R_{\pm}(s,w)=\frac{\epsilon_\pm a_{1-s}(2)}{2\cdot 4^w}\lz R(s,w;\psi_2,\psi_0)-2^w R(s,w;\psi_1,\psi_1)\pz,
\end{equation}
and by \eqref{eq:fe in w},
\begin{equation}\label{eq:R in terms of R+-}
	R(s)=\frac{Y_+\lz\tfrac12-\tfrac s3\pz}{3\zeta\lz1-\tfrac{2s}{3}\pz}R_+\lz\tfrac{2s}{3},\tfrac12+\tfrac s3\pz+\frac{Y_-\lz\tfrac12-\tfrac s3\pz}{3\zeta\lz1-\tfrac{2s}{3}\pz}R_-\lz \tfrac{2s}3,\tfrac 12+\tfrac s3\pz.
\end{equation}

Now we compute $R(s,w;\psi,\psi')$.
\begin{lemma}\label{le:R(s,w,psi,psi')}
	For $0<\re(s)<1$, $\re(w)>1/2$ and $\re(s+w)>1$, we have
	$$
\begin{aligned}
			R(s,w;\psi_1,\psi_1)&=\frac{\zeta(s)E(s)}{\zeta(2s)\zeta(2)}\lz1+\frac{2}{2^s(2-2^s)}\pz^{-1} P(s,w),\\
		R(s,w;\psi_2,\psi_0)&=R(s,w;\psi_1,\psi_1)(2^w-2^{-w})^{-1},
\end{aligned}
	$$ where
$$
\begin{aligned}
			E(s)&=\prod_p\lz 1+\frac{p^s-p^{1-s}+1}{(p+1)(p^s+1)(p^{1-s}-1)}\pz,\\
	P(s,w)&=\prod_{p>2}\lz1+\frac{p(p^{1+s}-p^{2s}+p)+p^{2w}(-p^{2+s}-p^2-p^{1+3s}+p^{1+s}+p^{4s})}{(p^{2s}-p^{1+s}-p)(p^{2w}-1)(p^{2s+2w}-p)}\pz.
\end{aligned}
$$
\end{lemma}
\begin{proof}
	By \eqref{eq:definition of B(s,w;t,psi,psi')} and Lemma \ref{le:continuation} with $N=1$, we have
	\begin{equation}\label{eq:expanded B}
\begin{aligned}
			B(s,w;t,\psi,\psi')&=\sum_{k\geq 1}\frac{\psi'(k)}{k^w}\frac{L(s,\leg{4k}{\cdot}\psi)L(s+t,\leg{4k}{\cdot}\psi)}{L(2s,\leg{4k}{\cdot}^2\psi^2)L(2s+2t,\leg{4k}{\cdot}^2\psi^2)}\\
			&\cdot\quad E_1\lz s,t;\leg{4k}{\cdot}\psi\pz P(s,t,k;\psi).
\end{aligned}
	\end{equation}
	If $\psi=\psi'=\psi_1$, we obtain the poles at $s+t=1$ from the summands with $k=\square$, and each such summand has a residue
	$$
		\begin{aligned}
			&\frac{\zeta_{(4k)}(s)\res{s+t=1}\zeta_{(4k)}(s+t)}{\zeta_{(4k)}(2s)\zeta_{(4k)}(2)}E_1\lz s,1-s;\leg{4k}{\cdot}\pz P(s,1-s,k;\psi_1)\\
			&=\frac{\zeta(s)}{\zeta(2s)\zeta(2)}E_1\lz s,1-s;\leg{4k}{\cdot}\pz P(s,1-s,k;\psi_1)\prod_{p|4k}\frac{\lz1-\frac{1}{p^s}\pz\lz1-\frac1p\pz}{\lz1-\frac{1}{p^{2s}}\pz\lz1-\frac{1}{p^2}\pz}\\
			&=\frac{\zeta(s)}{\zeta(2s)\zeta(2)}E_1\lz s,1-s;\leg{4k}{\cdot}\pz P(s,1-s,k;\psi_1)\prod_{p|4k}\lz1+\frac{1}{p^s}\pz^{-1}\lz1+\frac{1}{p}\pz^{-1}.
		\end{aligned}
	$$
Now using \eqref{eq:def of E}, we obtain for each $k=\square$
 $$
\begin{aligned}
	 	E_1\lz s,1-s;\leg{4k}{\cdot}\pz&=\prod_{p\nmid 4k}\lz 1+\frac{p^s-p^{1-s}+1}{(p+1)(p^s+1)(p^{1-s}-1)}\pz\\
	 	&=E(s)\prod_{p|4k}\lz 1+\frac{p^s-p^{1-s}+1}{(p+1)(p^s+1)(p^{1-s}-1)}\pz^{-1},
\end{aligned}
 $$
where we define
\begin{equation}\label{eq:E(s)}
	E(s)=\prod_p\lz 1+\frac{p^s-p^{1-s}+1}{(p+1)(p^s+1)(p^{1-s}-1)}\pz.
\end{equation}
Therefore
\begin{equation}\label{eq:expr R}
\begin{aligned}
		R(s,w;\psi_1,\psi_1)&=\frac{\zeta(s) E(s)}{\zeta(2s)\zeta(2)}\sum_{k\geq 1}\frac{\psi_1(k^2)P(s,1-s,k^2;\psi_1)}{k^{2w}}\\
		&\cdot\quad\prod_{p|4k}\lz1+\frac{1}{p^s}\pz^{-1}\lz1+\frac{1}{p}\pz^{-1}\lz1+\frac{p^s-p^{1-s}+1}{(p+1)(p^s+1)(p^{1-s}-1)}\pz^{-1}\\
		&=\frac{\zeta(s)E(s)}{\zeta(2s)\zeta(2)}\sum_{k\geq 1}\frac{\psi_1(k^2)P(s,1-s,k^2;\psi_1)}{k^{2w}}\prod_{p|4k}\lz1+\frac{p}{p^s(p-p^s)}\pz^{-1}.
\end{aligned}
\end{equation}
The sumands above contribute only if $k$ is odd, in which case \eqref{eq:def of P} and \eqref{eq:Gauss sum evaluation} yields
$$
\begin{aligned}
		&P(s,1-s,k^2;\psi_1)=\prod_{p|k}\lz 1+\lz1-\frac{1}{p^{1-s}}\pz^{-1}\sum_{j\geq 1}\frac{G\lz \bfrac{\cdot}{p^j},k^2\pz}{p^{js+j/2}}\pz\\
		&=\prod_{p^\alpha||k^2}\lz1+\lz1-\frac{1}{p^{1-s}}\pz^{-1}\sum_{j=1}^{\alpha/2}\frac{\phi(p^{2j})}{p^{2js+j}}+\lz1-\frac{1}{p^{1-s}}\pz^{-1}\frac{p^{\alpha+1/2}}{p^{(\alpha+1)s+\alpha/2+1/2}}\pz\\
		&=\prod_{p^\alpha||k^2}\lz1+\frac1{1-p^{s-1}}\sum_{j=1}^{\alpha/2}\frac{p^{2j-1}(p-1)}{p^{2js+j}} +\frac1{1-p^{s-1}}\cdot\frac{1}{p^{(\alpha+1)s-\alpha/2}}\pz\\
		&=\prod_{p^\alpha||k^2}\lz1+\frac1{1-p^{s-1}}\cdot\frac{p-1}{p}\cdot \frac{1}{p^{2s-1}}\sum_{j=0}^{\alpha/2-1}p^{j(1-2s)}+\frac1{1-p^{s-1}}\cdot\frac{1}{p^{(\alpha+1)s-\alpha/2}}\pz\\
		&=\prod_{p^\alpha||k^2}\lz1+\frac1{1-p^{s-1}}\cdot\frac{p-1}{p}\cdot \frac{1}{p^{2s-1}}\frac{p^{\alpha/2(1-2s)}-1}{p^{1-2s}-1}+\frac1{1-p^{s-1}}\cdot\frac{1}{p^{(\alpha+1)s-\alpha/2}}\pz\\
		&=\prod_{p^\alpha||k^2}\lz 1+\frac1{1-p^{s-1}}\cdot\frac{p-1}{p}\cdot\frac{p^{\alpha/2-\alpha s}-1}{1-p^{2s-1}}+\frac1{1-p^{s-1}}\cdot\frac{1}{p^{(\alpha+1)s-\alpha/2}}\pz\\
		&=\prod_{p^\alpha||k^2}\lz 1+\frac{p-1}{(p^s-p)(1-p^{2s-1})}+\frac{1}{p^{\alpha s-\alpha/2}}\lz\frac{1}{p^s(1-p^{s-1})}+\frac{p-1}{(p-p^s)(1-p^{2s-1})}\pz\pz.
\end{aligned}
$$
Therefore by expanding the sum over $k$ on the last line of \eqref{eq:expr R} into an Euler product and a straightforward computation (Mathematica can be handy), we find that it equals $\lz1+\frac{2}{2^s(2-2^s)}\pz^{-1}$ times
$$
\begin{aligned}
		&\prod_{p>2}\Bigg\{1+\lz1+\frac{p-1}{(p^s-p)(1-p^{2s-1})}\pz\lz1+\frac{p}{p^s(p-p^s)}\pz^{-1}\frac{1}{p^{2w}-1}\\
	&\quad\qquad+\frac{p+p^{1-s}}{p-p^{2s}}\lz1+\frac{p}{p^s(p-p^s)}\pz^{-1}\frac{1}{p^{2w+2s-1}-1}\Bigg\}\\
		&=\prod_{p>2}\lz1+\frac{p(p^{1+s}-p^{2s}+p)+p^{2w}(-p^{2+s}-p^2-p^{1+3s}+p^{1+s}+p^{4s})}{(p^{2s}-p^{1+s}-p)(p^{2w}-1)(p^{2s+2w}-p)}\pz
\end{aligned}
$$
Therefore
$$
	R(s,w;\psi_1,\psi_1)=\frac{\zeta(s)E(s)}{\zeta(2s)\zeta(2)}\lz1+\frac{2}{2^s(2-2^s)}\pz^{-1} P(s,w),
$$ where
$$
	P(s,w)=\prod_{p>2}\lz1+\frac{p(p^{1+s}-p^{2s}+p)+p^{2w}(-p^{2+s}-p^2-p^{1+3s}+p^{1+s}+p^{4s})}{(p^{2s}-p^{1+s}-p)(p^{2w}-1)(p^{2s+2w}-p)}\pz,
$$ and both $E(s)$ and $P(s,w)$ are absolutely convergent in our region.
\smallskip

If $\psi=\psi_2$ and $\psi'=\psi_0$, the residue in \eqref{eq:expanded B} comes from the summands $k=2\square.$ The computation up to \eqref{eq:E(s)} is completely analoguous (up to replacing the occurences of $\psi_1$ with $\psi_2$). For $P(s,1-s,2k^2;\psi_2)$, only the odd prime factors of $2k^2$ contribute. Therefore if $2k^2=2^{2a+1} k_0^2$ with $k_0$ odd, we have
$P(s,1-s,2k^2;\psi_2)=P(s,1-s,k_0^2;\psi_1)$, where we also used \eqref{eq:Gauss sum evaluation} to note that if $p$ is odd and $p^\alpha||k^2$, then $G\lz\leg{\cdot}{p^{\alpha+1}},2k^2\pz\psi_2(p^{\alpha+1})=G\lz\leg{\cdot}{p^{\alpha+1}},k_0^2\pz$).
Therefore
$$
	\begin{aligned}
		&R(s,w;\psi_2,\psi_0)=\frac{\zeta(s) E(s)}{\zeta(2s)\zeta(2)}\sum_{k\geq 1}\frac{P(s,1-s,2k^2;\psi_2)}{(2k^2)^w}\prod_{p|4k}\lz1+\frac{p}{p^{s}(p-p^s)}\pz^{-1}\\
		&=\frac{\zeta(s)E(s)}{\zeta(2s)\zeta(2)}\sum_{a\geq 0}\frac{1}{2^{(2a+1)w}}
		\sum_{k_0\geq 1}\frac{\psi_1(k_0^2)P(s,1-s,k_0^2;\psi_1)}{k_0^{2w}}\prod_{p|4k_0}\lz1+\frac{p}{p^{s}(p-p^s)}\pz^{-1}\\
		&=R(s,w;\psi_1,\psi_1)(2^w-2^{-w})^{-1}
	\end{aligned}
$$
\end{proof}

\begin{proof}[Proof of Proposition \ref{prop:residues} (3)]
From \eqref{eq:R+ as R(s,w,psi,psi')} and Lemma \ref{le:R(s,w,psi,psi')} we obtain
$$
	\begin{aligned}
		R_+(s,w)&=\frac{a_{1-s}(2)}{2\cdot 4^w} R(s,w;\psi_1,\psi_1)\lz (2^w-2^{-w})^{-1}-2^w\pz\\
		&=\frac{2^{s}(4^w-2)}{(2^{3w}-2^w)(2^{2s}-2^{s+1}-2)} \frac{\zeta(s)E(s)}{\zeta(2s)\zeta(2)}P(s,w),\\
		R_-(s,w)&=i R_+(s,w),
		\end{aligned}
		$$
	and \eqref{eq:R in terms of R+-} then gives
		$$
			\begin{aligned}
		R(s)&=\frac{Y_+\bfrac{3-2s}{6}}{3\zeta\lz1-\frac {2s}3\pz}\frac{2^{\frac{2s}3}(4^{\frac12+\frac s3}-2)}{(2^{\frac32+s}-2^{\frac12+\frac s3})(2^{\frac{4s}3}-2^{\frac{2s}3+1}-2)}\frac{\zeta\bfrac{2s}3E\bfrac{2s}3}{\zeta\bfrac{4s}3\zeta(2)}P\lz\frac{2s}{3},\frac{3+2s}{6}\pz\\
		&+\frac{iY_-\bfrac{3-2s}{6}}{3\zeta\lz1-\frac {2s}3\pz}\frac{2^{\frac{2s}3}(4^{\frac12+\frac s3}-2)}{(2^{\frac32+s}-2^{\frac12+\frac s3})(2^{\frac{4s}3}-2^{\frac{2s}3+1}-2)}\frac{\zeta\bfrac{2s}3E\bfrac{2s}3}{\zeta\bfrac{4s}3\zeta(2)}P\lz\frac{2s}{3},\frac{3+2s}{6}\pz\\
		&=\lz Y_+\bfrac{3-2s}{6}+iY_-\bfrac{3-2s}{6}\pz Y_+\bfrac{2s}{3} \frac{2^{\frac{2s}3}(4^{\frac12+\frac s3}-2)}{3(2^{\frac32+s}-2^{\frac12+\frac s3})(2^{\frac{4s}3}-2^{\frac{2s}3+1}-2)}\\
		&\quad \cdot\quad \frac{E\bfrac{2s}3}{\zeta\bfrac{4s}3\zeta(2)}P\lz\frac{2s}{3},\frac{3+2s}{6}\pz,
	\end{aligned}
$$ where we also used $\frac{\zeta(2s/3)}{\zeta(1-2s/3)}=Y_+(2s/3)$. 
We note that $P(2s/3,1/2+s/3)$ converges only for $\re(s)>1/2$. With the help of mathematica, we find that
$$
	\begin{aligned}
		P(s,w)&=\prod_{p>2}\lz1+\frac{p(p^{1+s}-p^{2s}+p)+p^{2w}(-p^{2+s}-p^2-p^{1+3s}+p^{1+s}+p^{4s})}{(p^{2s}-p^{1+s}-p)(p^{2w}-1)(p^{2s+2w}-p)}\pz\\
		&=\zeta_{(2)}(2s+2w-1)\prod_{p>2}\lz1-\frac{p^{1+2s}+p^{1+s}-p-p^{3s}}{p^s(-p^{1+s}-p+p^{2s})(p^{2w}-1)}\pz,
	\end{aligned}
$$
where the product is absolutely convergent at least for $\re(s)>0, \re(w)>1/2$. Hence for $\re(s)>0$, we obtain
$$
	P\lz\frac{2s}{3},\frac{3+2s}{6}\pz=\zeta_{(2)}(2s)\prod_{p>2}\lz1+\frac{p^{1+4s/3}+p^{1+2s/3}-p-p^{2s}}{p^{2s/3}(p^{1+2s/3}+p-p^{4s/3})(p^{1+4s/3}-1)}\pz.
$$

Therefore for $0<\re(s)<3/2$, $s\neq 1/2$, we have
$$
		R(s)=Y(s)\zeta(2s)Q_2(s)Q(s),
$$
where
$$
\begin{aligned}
	Y(s)&:=Y_+\bfrac{2s}{3}\lz Y_+\bfrac{3-2s}{6}+iY_-\bfrac{3-2s}{6}\pz,\\
	Q_2(s)&:=\lz1-\frac{1}{2^{2s}}\pz\frac{2^{\frac{2s}3}(4^{\frac12+\frac s3}-2)}{3(2^{\frac32+s}-2^{\frac12+\frac s3})(2^{\frac{4s}3}-2^{\frac{2s}3+1}-2)},\\
	Q(s)&:=\frac{E\bfrac{2s}{3}}{\zeta\bfrac{4s}{3}\zeta(2)}\prod_{p>2}\lz1+\frac{p^{1+4s/3}+p^{1+2s/3}-p-p^{2s}}{p^{2s/3}(p^{1+2s/3}+p-p^{4s/3})(p^{1+4s/3}-1)}\pz.
\end{aligned}
$$
\end{proof}
Theorem \ref{thm:MDS} now follows directly from Propositions \ref{prop:continuation of A final} and \ref{prop:residues}

\bibliography{../refs}
\bibliographystyle{abbrv}

\end{document}